arXiv:1301.5493v1 [math.DG] 23 Jan 2013# Yet More Smooth Mapping Spaces and Their Smoothly Local Properties

Andrew StaceyJanuary 23, 2013

**Abstract**

Motivated by the definition of the smooth manifold structure on a suitable mapping space, we consider the general problem of how to transfer local properties from a smooth space to an associated mapping space. This leads to the notion of *smoothly local* properties.

In realising the definition of a local property at a particular point it may be that there are choices that need to be made. To say that the local property is smoothly local is to say that those choices can be made smoothly dependent on the point. In particular, that a manifold has charts is a local property. A local addition is the structure needed to say that there is a way to choose a chart about each point so that it varies smoothly with that point.

To be able to extend these ideas beyond that of the local additivity of manifolds we work in a category of generalised smooth spaces. We are thus are able to consider more general mapping spaces than just those arising from the maps from one smooth manifold to another and thus able to generalise the standard result on when this space of maps is again a smooth manifold.

As applications of this generalisation we show that the mapping spaces involving the various figure 8s from String Topology are manifolds, and that they embed as submanifolds with tubular neighbourhoods in the corresponding loop spaces. We also show that applying the mapping space functor to a regular map of manifolds produces a regular map on the mapping spaces.
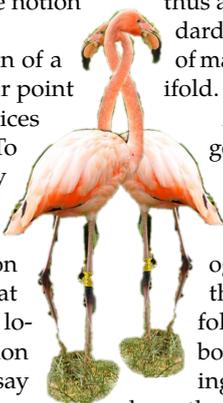

## Contents

**1 Introduction**     2

**2 Categorical Matters**     3

**3 Smoothly Local Properties**     8
   3.1 Structural Properties . . . . . . . . . . . . . . . . . . . . . . . . 8
   3.2 Deformation Properties . . . . . . . . . . . . . . . . . . . . . . 11
   3.3 Localisation Properties . . . . . . . . . . . . . . . . . . . . . . 15
   3.4 Miscellaneous Properties . . . . . . . . . . . . . . . . . . . . . 16
1

# 1 Introduction

> I'm smooth with a capital 'smoo'.
>
> The Cat, Red Dwarf

In this article we study the problem of transferring local information from a generalised smooth space to an associated mapping space. Specifically, if $X$ is a generalised smooth space with some local property and $S$ is another generalised smooth space, under what circumstances can we transfer that local property of $X$ to the mapping space $C^\infty(S, X)$?

The motivating example is the case where the target space is a smooth manifold. It is well-known that if $M$ and $N$ are finite dimensional smooth manifolds with $N$ compact then the space of smooth maps, $C^\infty(N, M)$, can be given the structure of an infinite dimensional smooth manifold modelled on Fréchet spaces (see [Mic80, KM97]). However, the charts on $C^\infty(N, M)$ do not come directly from the charts on $M$: there is a more involved process using something called a *local addition*. The difficulty with the charts on $M$ is that they only provide information about what is happening near one point. When considering a smooth map $\alpha \colon N \to M$ we need information simultaneously about $M$ near every point in the image of $\alpha$. A local addition can be viewed as a way of specifying charts at every point in $M$ so that they vary smoothly over $M$. This is the key information required to build the charts on $C^\infty(N, M)$.

This procedure generalises to the situation outlined in the first paragraph. If we want to transfer a local property from $X$ to $C^\infty(S, X)$ we need to show that $X$ can be given that local property in a way that varies smoothly over $X$. We choose to call this a *smoothly local* property. Such properties tend to behave well under taking mapping spaces in that there will be a relatively straightforward result giving simple conditions under which this smoothly local property on $X$ transfers to $C^\infty(S, X)$ from which it can be deduced that $C^\infty(S, X)$ has the original local property.

This does not mean that it is easy to show that if $X$ has a certain local property then $C^\infty(S, X)$ also has it. As we shall see in examples, the key step of showing that $X$ has the corresponding smoothly local property may involve some complicated pieces. However, this method does mean that the hard work is carried out on $X$ rather than on $C^\infty(S, X)$. In the motivating case of manifolds this means that the hard work is done in finite dimensions rather than infinite dimensions.

This paper is organised as follows. It breaks down into three main parts. In the first, Section 2, we shall explain the context in which we work. We have already used the term "generalised smooth space" and this needs clarifying and justifying. Also we need to explain how to make precise the notion of a *local* property for such an object. In the second part, Sections 3 and 4, we shall do the "easy" part of this paper wherein we define a variety of smoothly local properties and show how and when they transfer to mapping spaces. The final part, Sections 5 and 6, is concerned with applying the previous sections to particular situations, mostly that of finite dimensional manifolds. While the work of this section is not exactly hard, it uses more heavy-duty tools than the earlier ones.

Let us conclude this introduction by stating one of the theorems that we shall prove in this paper.

**Theorem 1.1** *Let M be a finite dimensional smooth manifold. Let S be a Frölicher space with the property that there is a non-zero smooth function $C^\infty(S, \mathbb{R}) \to \mathbb{R}$ with support in $C^\infty(S, (-1, 1))$. Then $C^\infty(S, M)$ is a smooth manifold which is locally modelled on its kinematic tangent spaces.*

*Suppose, in addition, that N is another finite dimensional smooth manifold and $f \colon M \to N$ a regular smooth map. Then $C^\infty(S, f) \colon C^\infty(S, M) \to C^\infty(S, N)$ is a regular smooth map.*

## 2 Categorical Matters

The original purpose of this article was to push the construction of the smooth manifold structure on mapping spaces even further than before. The standard construction begins with two finite dimensional smooth manifolds, say $M$ and $N$, with one of them, say $N$, compact. It then defines the structure of a Fréchet manifold on the set of smooth maps $C^\infty(N, M)$.

In the construction it soon becomes apparent that the majority of the structure on the source space, $N$ in the above, is not directly used. Its main role is to explain what are the smooth maps $N \to M$. Thus we can consider replacing $N$ by some other object so long as we still know what we mean by a smooth map $N \to M$. This suggests working in a category of generalised smooth spaces.

Doing so has implications for the rest of the work. By working in a category of generalised smooth spaces, the question "Is $X$ a smooth manifold?" can be broken down into "Is $X$ a smooth space?" and "Is that smooth space a smooth manifold?". The second is a question of *property* rather than *structure*. Moreover, the majority of known categories of generalised smooth spaces are cartesian closed whereupon for mapping spaces the first question can always be answered "Yes".

We therefore start by fixing a category of generalised smooth spaces. There are many of these and a survey of several is contained in [Sta11]. For this article the exact choice is immaterial.

Once we have chosen this category we need to know what it means to talk of "local" properties of an object. This suggests that we need some way to specify a topology on a smooth space. There is an important point to note when it comes to the relationship between topology and smoothness. Some of the categories define a smooth space as a topological space with some extra information. However, not all of them do so. Others define a smooth space as a set with some extra information. The distinction is as to whether differentiability is something that is built on top of continuity or alongside it. Whether or not our category came equipped with such we need to have a notion of topology to work with. Indeed, even if our category has a natural choice we need not use that one. We therefore need to choose not just a category of smooth spaces but also a functor to topological spaces.

We can make various assumptions about these choices. The initial part of this work does not need much. In fact, it is not necessary to be working with smooth spaces of any guise: almost any category will do. However, for the later stages we will need more structure. We shall assume this structure from the start.

Let us fix the data now. We choose the following things subject to the stated conditions.

1. A category, $\mathfrak{S}$, of smooth spaces.

    As mentioned above, initially we don't need anything more than the category. However, for later we shall want to assume various things about this category. These are that this category should fit the mould of [Sta11]. In short, it should be a topological concrete category. This implies that it is complete, cocomplete, cartesian closed, and concrete. The forgetful functor to Set is faithful and has a left and a right adjoint (the *discrete* and *indiscrete* smooth structures on a set, respectively).

    We assume that the category of smooth manifolds embeds fully and faithfully in $\mathfrak{S}$ (respecting the underlying sets) and we fix such an embedding. Via this embedding, we identify $\mathbb{R}$.

    We shall refer to objects in this category as *smooth spaces* and morphisms as *smooth maps*.

2. A functor, $\mathfrak{T} \colon \mathfrak{S} \to \mathrm{Top}$, to the category of topological spaces which respects underlying sets.

    We assume that $\mathfrak{T}(\mathbb{R})$ is the real line with its usual topology.

    We shall use $\mathfrak{T}$ to allow us to talk of topological matters in $\mathfrak{S}$, frequently without explicit reference to $\mathfrak{T}$.

The assumption that $\mathfrak{T}(\mathbb{R})$ is the usual topology on $\mathbb{R}$ puts limits on the possible topology functors. In particular, for a smooth space $X$, since every smooth map $\mathbb{R} \to X$, or smooth *curve* in $X$, must be continuous the topology on $\mathfrak{T}(X)$ can be no stronger than the strongest topology for which all the smooth curves are continuous. Conversely, since every smooth map $X \to \mathbb{R}$, or smooth *functional* on $X$, must be continuous the topology on $\mathfrak{T}(X)$ can be no weaker than the weakest topology for which all the smooth functionals are continuous. We therefore have two extremes.

**Definition 2.1** *The* curved topology functor $\mathfrak{T}_c \colon \mathfrak{S} \to \mathrm{Top}$ *is the functor which assigns to a smooth space X the strongest topology for which every smooth curve $\mathbb{R} \to X$ is continuous. We shall say that a subset in this topology is* curve open.

*The* functional topology functor $\mathfrak{T}_f \colon \mathfrak{S} \to \mathrm{Top}$ *is the functor which assigns to a smooth space X the weakest topology for which every smooth functional $X \to \mathbb{R}$ is continuous. We shall say that a subset in this topology is* functionally open.

*Given an arbitrary topology functor* $\mathfrak{T}\colon \mathcal{S} \to \mathrm{Top}$ *satisfying the assumption that* $\mathfrak{T}(\mathbb{R})$ *is the usual topology, we say that a smooth space X is* smoothly $\mathfrak{T}$–regular *if* $\mathfrak{T}(X) = \mathfrak{T}_f(X)$.

*We say that a smooth space is* smoothly regular *if it is smoothly $\mathfrak{T}_c$–regular.*

So for a smoothly regular space all topologies agree. Note that all finite dimensional smooth manifolds are smoothly regular, as are several important infinite dimensional ones.

The idea of a smoothly local property is that for a smooth space $X$ we are able to choose $\mathfrak{T}$–open neighbourhoods of each point with the given property in such a way that everything varies smoothly over $X$. Ignoring the property itself, we want to choose neighbourhoods of each point which vary smoothly over $X$. So for $x \in X$ we choose $V_x \subseteq X$ in a "smooth fashion".

That is a vague statement so let us make it more precise. To do so, let us consider the motivating example of smooth manifolds. For a smooth manifold, say $M$, the "smoothly local property" that we want is that it admit a *local addition*. This is a smooth map $\eta\colon TM \to M$ with certain properties. Key amongst those is that the map $\pi \times \eta\colon TM \to M \times M$ is a diffeomorphism onto a neighbourhood, say $V$, of the diagonal. Here, $\pi\colon TM \to M$ is the anchor map. When restricted to a fibre this means that $\eta_x\colon T_xM \to M$ is a diffeomorphism onto a neighbourhood of $x$, i.e., is a local chart at $x$. Another way to look at this is that there is a neighbourhood of the diagonal, $V \subseteq M \times M$, and a map $\tau = \pi_1\colon V \to M$ with the property that the fibres of this map can be given the structure of a vector space in such a way that this structure varies smoothly with the basepoint.

Thus if we ignore the vector space structure our starting point is a neighbourhood of the diagonal, $V \subseteq M \times M$, and the projection $\tau = \pi_1\colon V \to M$.

It does not take a large leap to generalise this a little further. The neighbourhood $V$ is a special case of a tubular neighbourhood. So it is likely that it would not take much effort to extend the concepts that we shall discuss from diagonals to submanifolds. The key difference is that in the case of a submanifold there is not a natural projection from the ambient manifold to the submanifold. But we only need this on the neighbourhood and we can stipulate this as a requirement.

Thus our definitions and results will all involve the following initial data in some form.

**Definition 2.2** *Let $X$ and $Y$ be smooth spaces with $X$ embedded in $Y$. A* smoothly local base *for $X$ in $Y$ is a pair $(V, \tau)$ where $V \subseteq Y$ is a subspace of $Y$ containing $X$ and $\tau\colon V \to X$ is a smooth map such that $V$ is a $\mathfrak{T}$–open neighbourhood of $X$ in $Y$ and $\tau$ is a retraction.*

*As shorthand, we shall write $X \stackrel{\frown}{\subseteq} V \subseteq Y$ for this structure.*

That is to say that when we apply the functor $\mathfrak{T}$ to everything in sight then $\mathfrak{T}(V)$ is open in $\mathfrak{T}(Y)$.

When we define a smoothly local property, starting with $X \stackrel{\frown}{\subseteq} V \subseteq Y$, our aim will be for this to be (relatively) easy to transfer to a mapping space. Let $S$ be a smooth space and consider the result of applying the exponential functor $C^\infty(S, -)$ to the above triple. We obtain smooth spaces $C^\infty(S, X)$, $C^\infty(S, V)$, $C^\infty(S, Y)$ and the embeddings and retraction on the original spaces define embeddings and a retraction on the mapping spaces. However, we do not necessarily get that $C^\infty(S, V)$ is $\mathfrak{T}$–open in $C^\infty(S, Y)$. A simple example (the topologies here are the "standard" ones) is where $X = \mathbb{R}$, $Y = \mathbb{R}^2$, and $V = \mathbb{R} \times (-1, 1)$. Then the topology on $C^\infty(\mathbb{R}, \mathbb{R}^2)$ is the compact–open topology (of arbitrary derivatives) which only gives control over what happens on a compact subset of $\mathbb{R}$. So it can only be guaranteed that a compact subset of $\mathbb{R}$ is mapped into $V$, not the whole of it. Thus $C^\infty(\mathbb{R}, V)$ is not open in $C^\infty(\mathbb{R}, \mathbb{R}^2)$.

As that example shows this is a matter more for the source space than the target. As our situation is extremely general at this point we can do no more than characterise and name the property that we want the source to have.

**Definition 2.3** *A smooth space, say S, is said to be* smoothly $\mathfrak{T}$–compact *if whenever $U \subseteq X$ is such that $\mathfrak{T}(U)$ is an open subset of $\mathfrak{T}(X)$ then $\mathfrak{T}(C^\infty(S, U))$ is an open subset of $\mathfrak{T}(C^\infty(S, X))$.*

*Let $\mathcal{C}$ denote the full subcategory of $\mathcal{S}$ consisting of smoothly $\mathfrak{T}$–compact spaces.*

This property shares many of the standard properties of compactness.

**Lemma 2.4** *The image of a smoothly $\mathfrak{T}$–compact space is smoothly $\mathfrak{T}$–compact.*

*Proof.* Let $S, T, X, U$ be smooth spaces. Suppose that $U \subseteq X$ is such that $\mathfrak{T}(U)$ is open in $\mathfrak{T}(X)$. Suppose that there is a smooth map $f\colon S \to T$ which is surjective. The map $f$ defines a map $f^*\colon C^\infty(T, X) \to C^\infty(S, X)$ and similarly for target $U$. As $f$ is surjective, if $g\colon T \to X$ is such that $g \circ f\colon S \to X$ has image in $U$ then $g$ must have image in $U$. Hence $C^\infty(T, U) = f^{*-1} C^\infty(S, U)$ and thus is open in $C^\infty(T, X)$.

In Section 6 we shall see what this looks like when the category of smooth spaces is that of Frölicher spaces and the topology functor is the curved topology functor. As we shall see, in that situation a compact manifold is smoothly $\mathfrak{T}$–compact.

Then "by definition" we have the next result.

**Lemma 2.5** *Let $X \overset{\frown}{\subseteq} V \subseteq Y$ be a smoothly local base for $X$ in $Y$. Let $S$ be a smoothly $\mathfrak{T}$–compact space. Then $C^\infty(S, X) \overset{\frown}{\subseteq} C^\infty(S, V) \subseteq C^\infty(S, Y)$ is a smoothly local base for $C^\infty(S, X)$ in $C^\infty(S, Y)$.*

With the functional topology we can introduce a slightly weaker concept than $\mathfrak{T}_f$–compactness. A subspace $U$ of a smooth space $X$ is $\mathfrak{T}_f$–open if for every $p \in U$ there is some smooth $f_p\colon X \to [0, 1]$ with $f_p(p) = 1$ and $f_p^{-1}((0, 1]) \subseteq U$. In dealing with mapping spaces our common theme is that when a choice can be made it is best if it is possible to make it smoothly.

**Definition 2.6** *Let $X$ be a smooth space. A subspace $U \subseteq X$ is said to be* smoothly functionally open *if there is a smooth map $f\colon U \to C^\infty(X, [0, 1])$ with the property that for each $p \in U$ then $f(p)(p) = 1$ and $f(p)^{-1}((0, 1]) \subseteq U$.*

*A smooth space $S$ is said to be* weakly $\mathfrak{T}_f$–compact *if whenever $U \subseteq X$ is smoothly functionally open then $C^\infty(S, U) \subseteq C^\infty(S, X)$ is smoothly functionally open.*

Equivalently, we could relax this to allow codomain $\mathbb{R}$ or $f(p)(p) \neq 0$. The set of smoothly functionally open subspaces does not (necessarily) form a topology: it may not be closed under arbitrary unions. The nearest topology is the functional topology.

**Lemma 2.7** *A smooth space $S$ is weakly smoothly compact if and only if there is a smooth function $g\colon C^\infty(S, \mathbb{R}) \to [0, 1]$ with the property that $g(0) = 1$ and $g^{-1}((0, 1]) \subseteq C^\infty(S, (-1, 1))$.*

*Proof.* The "only if" follows from the fact that $(-1, 1)$ is smoothly functionally open in $\mathbb{R}$, whence by definition $C^\infty(S, (-1, 1))$ is smoothly functionally open in $C^\infty(S, \mathbb{R})$ and thus is open for the functional topology, whence there is a smooth function as stated.

For the "if" part, let $U \subseteq X$ be smoothly functionally open. Let $f \colon U \to C^\infty(X, [0, 1])$ be a suitable smooth function. Via adjunction we consider this a smooth function $f \colon U \times X \to [0, 1]$. Taking mapping spaces, we get a smooth function $C^\infty(S, U \times X) \to C^\infty(S, [0, 1])$ and thus $F \colon C^\infty(S, U) \times C^\infty(S, X) \to C^\infty(S, [0, 1])$. Via the inclusion $[0, 1] \to \mathbb{R}$ we regard this as a map into $C^\infty(S, \mathbb{R})$ which we then compose with the function $\beta \mapsto g(1 - \beta)$ to get $G \colon C^\infty(S, U) \times C^\infty(S, X) \to [0, 1]$.

Let $\alpha \colon S \to U$ and $\beta \colon S \to X$ be smooth maps. Then $F(\alpha, \beta)$ is the map $s \mapsto f(\alpha(s), \beta(s))$. So if $\beta = \alpha$ then $F(\alpha, \beta)$ is the constant map at 1 and if there is some $s \in S$ for which $\beta(s) \notin U$ then for that $s$, $F(\alpha)(\beta)(s) = 0$. Thus if $\beta = \alpha$ then $1 - F(\alpha, \beta)$ is the zero map and so $G(\alpha, \beta) = 1$. Whilst if $\beta$ strays outside $U$ then $1 - F(\alpha, \beta)$ takes on the value 1, whence $G(\alpha, \beta) = 0$. Hence $G$ shows that $C^\infty(S, U)$ is smoothly functionally open in $C^\infty(S, X)$ as required. 🐦

Note that the existence of such a function, say $g$, on $C^\infty(S, \mathbb{R})$ implies that every smooth function on $S$ has compact image. This is because for a smooth function $f \colon S \to \mathbb{R}$ then the map $\mathbb{R} \times S \to \mathbb{R}$ given by $(t, s) \mapsto tf(s)$ is again smooth, whence $t \mapsto tf$ is a smooth path, say $\alpha$, in $C^\infty(S, \mathbb{R})$. At 0 it is the constant function 0 and so $g \circ \alpha(0) = 1$. There is, therefore, some $t \neq 0$ such that $g \circ \alpha(t) \neq 0$ whence $tf \in C^\infty(S, (-1, 1))$ and hence, as $t \neq 0$, $f$ is bounded. As this holds for all $f$, we must have that each $f$ has compact image.

Let us end this section with a comment on nomenclature. For a particular property that we are interested in we could provide several variants of the basic definition. Of these, we have tried to select the most practical for the most straightforward name.

The first variation is because most of our properties involve a choice of some structure. Thus we can consider spaces which admit a choice and spaces together with a specific choice. Although both have uses, the former is the one we anticipate being used most.

The second variation is from the nature of local properties. In topology, when a space is said to be "locally $P$" then the usual meaning is that for every point, every neighbourhood of that point contains a neighbourhood with property $P$. In differential topology, the purpose of the local regions is that it is nicer to work there than over the whole space. It is sometimes enough, therefore, to simply have some nice neighbourhood. Moreover, when dealing with mapping spaces it is not true that every open set in the mapping space is formed by taking an open set in the original space and exponentiating.

We shall use the template "$X$ is smoothly locally $P$ (in $Y$)" where if the "in $Y$" is omitted then it is understood that $Y = X \times X$ and $X$ is the diagonal (technically, $X$ is identified with its image under the diagonal map but we shall rarely make that distinction). What we will mean by this will depend on what other information has been given. In general, it will mean that there is a neighbourhood, say $V$, of $X$ in $Y$ with a retraction $V \to X$ such that $X \mathrel{\widehat{\subseteq}} V \subseteq Y$ has some property.

The last variation on the definition is as to the amount of information that has to be chosen. It may be that either $V$ or the retraction has been given beforehand, in which case we shall assume that those are the ones in use. That is, we trust to context to distinguish these situations from those where no information is given.

In particular, if $Y$ is not given (so is assumed to be $X \times X$) then the retraction will always be assumed to be the restriction of the projection $\pi_1 \colon X \times X \to X$ onto the first factor.

## 3 Smoothly Local Properties

In this section we shall define some properties that we wish to study. Recall that our properties will start with a triple $X \mathrel{\widehat{\subseteq}} V \subseteq Y$ of smooth spaces where $V$ is an open neighbourhood of $X$ in $Y$. The properties can be grouped together according to type.

The first type, structural properties, are to do with the existence of structure in terms of the existence of morphisms satisfying properties. Given a suitable pair of smooth spaces, $X$ and $Y$, and a smoothly local neighbourhood base, $X \mathrel{\widehat{\subseteq}} V \subseteq Y$ with retraction $\tau \colon V \to X$, a structural property is one on the fibres of $\tau$. The key example of this is that of being smoothly locally additive wherein the fibres of $\tau$ become vector spaces. In the example of $Y = X \times X$ these provide the charts for a smooth manifold structure. These are primarily of use when considering the structure of mapping spaces with a single, fixed source. Thus we start with $X \subseteq Y$ and look for information about $C^\infty(S, X) \subseteq C^\infty(S, Y)$.

The second type, deformation properties, are to do with moving $X$ in $Y$ along the fibres of $\tau$. The idea here is to consider maps $\alpha \colon S \to Y$ which take some suitable subspace, say $T$, into $V$ and deform them so that $T$ is actually mapped into $X$. The deformations involved have to be globally defined on $Y$ but only affect the region in $V$. These properties allow us to consider what happens when we vary the source. Thus as well as $X \subseteq Y$ we also have two source spaces, say $S$ and $T$, which are related in some fashion. Then we consider the relationship of the relative mapping space, $C^\infty((S, T), (Y, X))$, to $C^\infty(S, Y)$. A key example of this is the deformation in string topology where a smooth map of two disjoint circles is deformed so as to bring the basepoints together.

The third type, localisation properties, is concerned with finding a function from $Y$ to $\mathbb{R}$ that "detects" $X$ and its smoothly local base. That is, finding a function supported in $V$ and value 1 on $X$. As this is concerned with maps from $Y$ the existence of such for $Y$ is less closely related to that on the mapping space.

Lastly, we consider how fibre bundles fit into this scheme. This is not so closely related to the other properties but as we shall see when we consider manifolds it is something that occurs as a useful tool at many stages.

### 3.1 Structural Properties

The structural properties all depend on the existence of certain morphisms satisfying certain identities. That they are preserved by the mapping space construction is fairly straightforward. Although our primary focus is on spaces which *admit* a particular structure we shall have use for spaces with a choice of such structure and so we start with these.

There are obviously many properties that we could define here. We pick three. The third is the key property relating to manifolds. The second will show its importance when we come to talk of tangent bundles. The first is closely related to the second and is more common in use. Obviously, the third implies the second which implies the first, and all depend on the fact that we know what the real numbers are in our category (whence, if armed with that knowledge, it would be simple to adapt these definitions to fairly arbitrary categories).

**Definition 3.1** *Let X and Y be smooth spaces with X embedded into Y.*

1. A smoothly local contraction *of X in Y consists of a smoothly local base, $X \stackrel{\frown}{\subseteq} V \subseteq Y$, and a smooth function $H \colon [0,1] \times V \to V$ with the following properties (writing $H_t(y)$ for $H(t,y)$ and where $\tau \colon V \to X$ is the retraction):*

   (a) $\tau H_t = \tau$,
   
   (b) $H_1 = Id$,
   
   (c) $H_0 = \tau$.

   *A morphism of smoothly local contractions $(X_1, Y_1, V_1, H_1) \to (X_2, Y_2, V_2, H_2)$ is a smooth map $g \colon Y_1 \to Y_2$ which intertwines the structures. This defines a category of pairs of spaces with a smoothly local contraction.*

   *We say that X is* smoothly locally contractible *in Y if it admits a smoothly local contraction.*

2. A smoothly local deflation *of X in Y consists of a smoothly local base, $X \stackrel{\frown}{\subseteq} V \subseteq Y$, and a smooth function $H \colon \mathbb{R} \times V \to V$ with the following properties (again writing $H_t(y)$ for $H(t,y)$):*

   (a) $\tau H_t = \tau$,
   
   (b) $H_s(H_t(y)) = H_{st}(y)$,
   
   (c) $H_1 = Id$,
   
   (d) $H_0 = \tau$

   *A morphism of smoothly local deflations $(X_1, Y_1, V_1, H_1) \to (X_2, Y_2, V_2, H_2)$ is a smooth map $g \colon Y_1 \to Y_2$ which intertwines the structures. This defines a category of pairs of spaces with a smoothly local deflation.*

   *We say that X is* smoothly locally deflatable *in Y if it admits a smoothly local deflation.*

3. A smoothly local addition *of X in Y consists of a smoothly local base, $X \stackrel{\frown}{\subseteq} V \subseteq Y$, together with smooth functions $\lambda \colon \mathbb{R} \times V \to V$ and $\alpha \colon V \times_X V \to V$ which together put the structure of a vector space on the fibres of $\tau \colon V \to X$ (with the inclusion $X \to V$ as the fibrewise zero vector).*

   *A morphism of smoothly local additions $(X_1, Y_1, V_1, \lambda_1, \alpha_1) \to (X_2, Y_2, V_2, \lambda_2, \alpha_2)$ is a smooth map $g \colon Y_1 \to Y_2$ which intertwines the structures. This defines a category of pairs of spaces with a smoothly local addition.*

   *We say that X is* smoothly locally additive *in Y if it admits a smoothly local addition.*

So in each case we have a category of pairs of smooth spaces with that structure and an obvious forgetful functor to the category of pairs of smooth spaces. The pairs in the image of that forgetful functor are the ones we are particularly interested in.

As all of these properties are of the form "exist maps satisfying identities", they are inherited by mapping spaces.

**Theorem 3.2** *Recall that $C$ is the category of smoothly $\mathfrak{T}$–compact spaces. Let $\mathcal{P}$ be the category corresponding to one of the structures in Definition 3.1. Then the exponential construction on $\mathcal{S}$ lifts to a functor:*

$$C \times \mathcal{P} \to \mathcal{P}$$

In particular, this applies to the property of being smoothly locally additive. This is very closely related to being a manifold. As we are working in a category of generalised smooth spaces, we need to define "being a manifold" as a *property* rather than a *structure*.

We start with the notion of a smooth vector space. Recall that we have identified $\mathbb{R}$ in our category of smooth spaces via the embedding from the category of finite dimensional smooth manifolds. This means that it is a ring object, whereupon we can define the notion of a module over it.

**Definition 3.3** *A* smooth vector space *is a smooth space, say $W$, together with smooth maps $\alpha \colon W \times W \to W$, $\lambda \colon \mathbb{R} \times W \to W$, and $0 \colon \{\star\} \to W$ which satisfy the identities for $W$ to be a module over $\mathbb{R}$.*

With this, we can define what it means for a smooth space to be a smooth manifold. As the definition requires a topology we use the functor $\mathfrak{T}$ to supply it. This influences our choice of name.

**Definition 3.4** *A $\mathfrak{T}$–pre-manifold in $\mathcal{S}$ is a smooth space $X$ with the property that each point $x \in X$ admits a $\mathfrak{T}$–open neighbourhood which is diffeomorphic to an open subset of a smooth vector space.*

*A $\mathfrak{T}$–manifold in $\mathcal{S}$ is a $\mathfrak{T}$–pre-manifold which is $\mathfrak{T}$–Hausdorff and $\mathfrak{T}$–paracompact.*

The precise relation between this definition and the usual one will depend somewhat on the category of smooth spaces and the topological functor in use. In all known cases, an ordinary smooth manifold will be a $\mathfrak{T}$–manifold but there may be smooth spaces which are $\mathfrak{T}$–manifolds but aren't "ordinary" manifolds.

As we are working in a category of smooth spaces it would be nice to replace the two purely topological conditions by corresponding smooth ones. A natural one for Hausdorff is *smoothly Hausdorff* meaning that smooth functions separate points. The main use of paracompactness is to transfer structure from local to global. This, we contend, is something that smoothly locally additive could conceivably do.

With the above definition we can deduce the following relationship between smoothly locally additive and being a manifold.

**Proposition 3.5** *Let $X$ be a smoothly locally additive space. Then $X$ is a $\mathfrak{T}$–pre-manifold.*

*Proof.* Recall that as we do not have an embedding space we assume it to be $X \times X$. So to say that $X$ is smoothly locally additive is to say that there is a neighbourhood, say $V$, of the diagonal in $X \times X$ with the property that each fibre of $\pi_1 \colon V \to X$ can be given the structure of a smooth vector space.

Let $x \in X$ and consider the map $i_x \colon X \to X \times X$ defined by $i_x(x') = (x, x')$. Then $V_x := i_x^{-1}(V)$ is an open set in $X$ and, by construction, contains $x$. The induced map $V_x \to V$ is (by construction) an embedding and its image is precisely $\pi_1^{-1}(x)$. We can therefore transfer the vector space structure to $V_x$ and thus $x$ has a neighbourhood that is diffeomorphic to a smooth vector space.

**Corollary 3.6** *If $X$ is smoothly locally additive then so is $C^\infty(S, X)$, whence $C^\infty(S, X)$ is a $\mathfrak{T}$–pre-manifold.*

There is a particularly interesting case of a morphism of smoothly locally additive spaces that we would like to introduce here.

**Definition 3.7** *Let $X$ and $Y$ be smooth spaces and $f \colon X \to Y$ a smooth map. We say that $f$ is* smoothly locally regular *if it fulfils the following conditions.*

1. *That $X$ is smoothly locally additive.*

2. *That the graph of $f$ is smoothly locally additive in $X \times Y$, where the retraction is the map $(x, y) \mapsto (x, f(x))$.*

3. *That there are choices of smoothly locally additive structures on these spaces such that the map $1 \times f \colon X \times X \to X \times Y$ preserves these structures.*

An immediate corollary of Theorem 3.2 is that this is preserved under taking mapping spaces.

**Corollary 3.8** *If $f \colon X \to Y$ is smoothly locally regular and $S$ is a smoothly $\mathfrak{T}$–compact space then $C^\infty(S, f) \colon C^\infty(S, X) \to C^\infty(S, Y)$ is smoothly locally regular.*

In Section 5 we shall relate this to the usual definition of a regular map between smooth manifolds.

## 3.2 Deformation Properties

In this section we consider properties to do with deforming $X$ within $Y$. Given a smoothly local base $X \mathrel{\widehat{\subseteq}} V \subseteq Y$ with retraction $\tau \colon V \to X$ we can consider deforming $V$ within $Y$ in such a manner that points move along the fibres of $\tau$. With just the $(X, V, Y)$ and the map $\tau \colon V \to X$ then any map $\alpha \colon S \to V$ can be projected to $X$ via $\alpha \mapsto \tau \circ \alpha$. But suppose we have a subset $T \subseteq S$ and a map $\alpha \colon S \to Y$ for which $\alpha(T) \subseteq V$. Then we might want to deform $\alpha$ such that $\alpha(T) \subseteq X$. Simply applying $\tau$ doesn't work because that can only be applied to $\alpha|_T$. What we need to do is deform $Y$ in such a way that $\alpha|_T$ ends up as $\tau \circ \alpha$. This deformation then tells us what to do with the rest of $\alpha$.

For the constructions in this section to work we need to make some assumptions on how smooth maps can be constructed. We will want to construct a smooth map by gluing together two smooth maps in a sheaf-like manner. We can either assume at the outset that our category of smooth spaces satisfies this property or work with only those smooth spaces where it is true. We will only need the sheaf-like property where the open sets are defined by smooth functionals.

Let us proceed with the definition of smoothly local deformability.

**Definition 3.9** *Let $X \subseteq Y$ be smooth spaces. A* smoothly local deformation *of $X$ in $Y$ consists of a smoothly local base, $X \mathrel{\widehat{\subseteq}} V \subseteq Y$ with retraction $\tau$, and a smooth map $\psi \colon \mathbb{R} \times V \to \mathrm{Diff}(Y)$ with the following properties, writing $\psi_{t,v}$ for $\psi(t, v)$.*

1. For $v \in V$, $t \mapsto \psi_{t,v}$ is a group homomorphism $\mathbb{R} \to \mathrm{Diff}(Y)$ with $\mathbb{R}$ additive.

2. For $t \in \mathbb{R}$ and $v \in V$, $\psi_{t,v}$ is the identity outside $V$.

3. For $t \in \mathbb{R}$, $v \in V$, and $y \in V$, $\tau\psi_{t,v}(y) = \tau(y)$.

4. For $v \in V$, $\psi_{1,v}(v) = \tau(v)$.

*We say that $X$ is* smoothly locally deformable *in $Y$ if it admits a smoothly local deformation in $Y$.*

Thus $\psi_{t,v}$ deforms $Y$ but only actually affects $V$. Within $V$, it only moves points along the fibres of $\tau\colon V \to X$. Lastly, $\psi_{1,v}$ is such that $v$ is taken to $\tau(v)$.

We have not defined a category of smoothly locally deformations as the interaction with morphisms is considerably more complicated and much less useful due to the use of the diffeomorphism group.

As mentioned above the real use of this is for considering maps $S \to Y$ where only a subset is required to lie in $X$.

**Proposition 3.10** *We consider the following ingredients: smooth spaces $X$, $Y$, $T$, $S$ with $X \subseteq Y$ and $T \subseteq S$; and a smooth map $\sigma\colon S \to [0,1]$. These must satisfy the following properties.*

1. *There is a smoothly local deformation of $X$ in $Y$, say with smoothly local base $X \overset{\circ}{\subseteq} V \subseteq Y$ and retraction $\tau\colon V \to X$.*

2. *$S$ is smoothly $\mathfrak{T}$–compact.*

3. *$\sigma(T) \subseteq \{1\}$.*

4. *There is a smoothly local base of $T$ in $S$, say $T \overset{\circ}{\subseteq} U \subseteq S$ with retraction $\nu\colon U \to T$, such that $\overline{\sigma^{-1}(0,1]} \subseteq U$.*

5. *Let $W := S \setminus \overline{\sigma^{-1}(0,1]}$. Then there is a pull-back square:*

$$\begin{array}{ccc} C^\infty(S,Y) & \longrightarrow & C^\infty(U,Y) \\ \downarrow & \lrcorner & \downarrow \\ C^\infty(W,Y) & \longrightarrow & C^\infty(U \cap W, Y) \end{array}$$

*Let $C^\infty((T,S),(X,Y))$ be the relative mapping space of smooth maps $S \to Y$ which take $T$ into $X$ and let $C^\infty((T,S),(V,Y))$ be defined similarly. Under the above conditions there is a diffeomorphism:*

$$C^\infty((T,S),(V,Y)) \to C^\infty((T,S),(X,Y)) \times_{C^\infty(T,X)} C^\infty(T,V)$$

*where $C^\infty((T,S),(X,Y)) \to C^\infty(T,X)$ is the restriction map and $C^\infty(T,V) \to C^\infty(T,X)$ is $C^\infty(T,\tau)$.*

*Proof.* We assume the notation of the question.

As we shall be discussing maps between mapping spaces and what happens to particular elements — i.e., maps — we shall use the term *map* for when we are talking about the maps between the mapping spaces and *function* for the elements thereof.

We define maps $C^\infty((T,S),(V,Y)) \to C^\infty(U,Y)$ and $C^\infty((T,S),(V,Y)) \to C^\infty(W,Y)$ as follows. For a smooth function $\alpha \colon S \to Y$ such that $\alpha(T) \subseteq V$, we define:

$$\alpha_U \colon U \to Y$$
$$s \mapsto \psi_{\sigma(s),\alpha(\nu(s))}(\alpha(s)),$$
$$\alpha_W \colon W \to Y$$
$$s \mapsto \alpha(s)$$

For $s \in U \cap W$ then $\sigma(s) = 0$ whence $\alpha_U(s) = \psi_{0,\alpha(\nu(s))}(\alpha(s)) = \alpha_W(s)$. Hence, by the pull-back assumption, these maps define a unique map $C^\infty((T,S),(V,Y)) \to C^\infty(S,Y)$. By construction, a function in the image of this map takes $T$ into $X$ and so the image of this map is in $C^\infty((T,S),(X,Y))$. The map $C^\infty((T,S),(V,Y)) \to C^\infty(T,V)$ is defined by restriction. We therefore have a map:

$$C^\infty((T,S),(V,Y)) \to C^\infty((T,S),(X,Y)) \times C^\infty(T,V).$$

With $\alpha$ and $\alpha_U$ as above then for $s \in T$, $\alpha_U(s) = \psi_{\sigma(s),\alpha(\nu(s))}(\alpha(s)) = \psi_{1,\alpha(s)}(\alpha(s)) = \tau(s)$. Hence the image of the above map is the subspace $\{(\beta,\gamma) \colon \beta|_T = \tau \circ \gamma\}$ and thus we have a map:

$$C^\infty((T,S),(V,Y)) \to C^\infty((T,S),(X,Y)) \times_{C^\infty(T,X)} C^\infty(T,V).$$

The inverse map is obtained in a similar fashion. Given $\beta \in C^\infty((T,S),(X,Y))$ and $\gamma \in C^\infty(T,V)$ such that $\beta|_T = \tau \circ \gamma$, we define:

$$\alpha_U \colon U \to Y$$
$$s \mapsto \psi_{-\sigma(s),\gamma(\nu(s))}(\beta(s)),$$
$$\alpha_V \colon W \to Y$$
$$s \mapsto \beta(s).$$

Exactly as before we construct out of this a smooth map $C^\infty((T,S),(X,Y)) \times_{C^\infty(T,X)} C^\infty(T,V) \to C^\infty(S,Y)$ which happens to lie in $C^\infty((T,S),(V,Y))$. It is inverse to the other map by construction. ⌣

**Corollary 3.11** *We assume the conditions of Proposition 3.10. Suppose that $X$ has one of the smoothly local structural properties of Section 3.1 in $Y$ where the neighbourhood and retraction are the same as that used in the construction in Proposition 3.10. Then $C^\infty((T,S),(X,Y))$ inherits this property in $C^\infty(S,Y)$.* ⌣

This deals with the case where $T \to S$ is an embedding. It is also possible to use smoothly local deformability to study the case where $T \to S$ is a quotient mapping. In this case it is of little use having a pair, say $X \subseteq Y$, on the target side since if $T \to S$ is a quotient and $\alpha \colon S \to Y$ is such that $\alpha$ takes $T$ to $X$ then $\alpha$ must take $S$ to $X$ as well. So we are interested in the relationship between $C^\infty(S,X)$ and $C^\infty(T,X)$.

As we no longer have an ambient space in which $X$ lies we consider it as a subspace of $X \times X$. Saying that $X$ is smoothly locally deformable means that when we have a pair $(x,y) \in X \times X$ that are sufficiently near the diagonal we can deform $X \times X$ to bring $(x,y)$ onto the diagonal. In this situation our projection is that onto the first factor (followed by the diagonal inclusion) so the deformation takes $(x,y)$ to $(x,x)$. Moreover, the deformation is along the fibres of the projection and so defines a deformation of $X$ itself. That is, with the notation of 3.9, since $\tau \colon X \times X \to X$ is globally defined, properties 2 and 3 mean that for $t \in \mathbb{R}$, $v \in V$, and any $x, y \in X$, $\tau \psi_{t,v}(x,y) = \tau(x,y) = x$. Thus we can use $\psi$ to define a map $\varphi \colon \mathbb{R} \times V \to \mathrm{Diff}(X)$. The properties of $\psi$ translate to the following properties of $\varphi$:

1. For $v \in V$, $t \mapsto \varphi_{t,v}$ is a group homomorphism.

2. For $t \in \mathbb{R}$ and $v = (x,y) \in V$, $\varphi_{t,v}$ is the identity outside $V_x = \{x' \in X : (x,x') \in V\}$.

3. For $v = (x,y) \in V$, $\varphi_{1,v}(y) = x$.

Using this reformulation, we can analyse the situation where $S$ is formed from $T$ by "bringing points together".

**Proposition 3.12** *We consider the following ingredients: smooth spaces $X$, $R$, and $T$ with $R \subseteq T$; and smooth maps $\rho \colon R \to T$ (in addition to the inclusion) and $\sigma \colon T \to [0, 1]$. These must satisfy the following properties.*

1. *$X$ is smoothly locally deformable. Let $V \subseteq X \times X$ be the corresponding neighbourhood of the diagonal.*

2. *$T$ is smoothly $\mathfrak{T}$–compact.*

3. *There is a smoothly local base for $R$ in $T$, say $R \overset{\frown}{\subseteq} U \subseteq T$ with retraction $\tau \colon U \to R$, such that $\overline{\sigma^{-1}(0,1]} \subseteq U$.*

4. *$\sigma(R) \subseteq \{1\}$.*

5. *The image of $\rho$ is disjoint from $\overline{\sigma^{-1}(0,1]}$.*

6. *Let $W := T \setminus \overline{\sigma^{-1}(0,1]}$. Then there is a pull-back square:*

$$\begin{array}{ccc} C^\infty(T, Y) & \longrightarrow & C^\infty(U, Y) \\ \downarrow & \lrcorner & \downarrow \\ C^\infty(W, Y) & \longrightarrow & C^\infty(U \cap W, Y) \end{array}$$

*Let $S$ be the coequaliser of $R \rightrightarrows T$ where the two maps are the inclusion and $\rho$. Let $C^\infty((R,T),(V,X))$ be the space of smooth maps $\alpha \colon T \to X$ with the property that for $r \in R$, $(\alpha(\rho(r)), \alpha(r)) \in V$. Then there is a diffeomorphism*

$$C^\infty((R,T),(V,X)) \to C^\infty(S, X) \times_{C^\infty(R,X)} C^\infty(R, V).$$

*Proof.* We assume the notation of the statement. Let $\varphi \colon \mathbb{R} \times V \to \mathrm{Diff}(X)$ be the map derived from the definition of smoothly locally deformable as discussed before the statement of the proposition. We define maps $C^\infty((R,T),(V,X)) \to$ $C^\infty(U, X)$ and $C^\infty((R, T), (V, X)) \to C^\infty(W, X)$ as follows. For a smooth function $\alpha \colon T \to X$ such that $(\alpha \circ \rho, \alpha)(R) \subseteq V$, we define:

$$\alpha_U \colon U \to X$$
$$t \mapsto \varphi_{\sigma(t), (\alpha \circ \rho \circ \tau(t), \alpha \circ \tau(t))}(\alpha(t)),$$
$$\alpha_W \colon W \to X$$
$$t \mapsto \alpha(t).$$

For $t \in U \cap W$ then $\sigma(t) = 0$ whence $\alpha_U(t) = \alpha(t)$. Hence, by the pull-back assumption, these maps define a unique map $C^\infty((R, T), (V, X)) \to C^\infty(T, X)$. As $\rho(R) \subseteq W$ and $R \subseteq U$, for $r \in R$ we have that $\alpha_U(r) = \alpha_W \circ \rho(r)$ and thus the image of the map lies in $C^\infty(S, X)$. We also get a map to $C^\infty(R, V)$ by sending $\alpha$ to $(\alpha \circ \rho, \alpha)$. Put together this yields:

$$C^\infty((R, T), (V, X)) \to C^\infty(S, X) \times C^\infty(R, V).$$

Both factors in the product have maps to $C^\infty(R, X)$. The first by composition with the quotient map $T \to S$ followed by restriction to $R$, the second by projection $V \to X$ onto the first factor. In both cases the result is $\alpha \circ \rho$, whence the above map is actually a map

$$C^\infty((R, T), (V, X)) \to C^\infty(S, X) \times_{C^\infty(R, X)} C^\infty(R, V).$$

Conversely, suppose that we have functions $\beta \colon S \to X$ and $\gamma \colon R \to V$ which agree when mapped to $C^\infty(R, X)$, we define:

$$\alpha_U \colon U \to X$$
$$t \mapsto \varphi_{-\sigma(t), \gamma \circ \tau(t)}(\beta(t)),$$
$$\alpha_W \colon W \to X$$
$$t \mapsto \beta(t).$$

Exactly as before this defines a smooth map $T \to X$. By construction it lies in $C^\infty((R, T), (V, X))$. Also by construction it is inverse to the above map. 🙂

## 3.3 Localisation Properties

In this section we consider properties to do with restriction and extension between $X$ and $Y$ and how these are inherited by mapping spaces. The primal situation is that of functionals and for those the basic tool is that of bump functions. Obviously the best situation to be in is to have plenty of bump functions.

**Definition 3.13** *Let $X$ and $Y$ be smooth spaces with $X \subseteq Y$. We shall say that a smoothly local base of $X$ in $Y$, say $X \hat{\subseteq} V \subseteq Y$, is* regular *if there is a smooth function $f \colon Y \to \mathbb{R}$ which is 1 on $X$ and supported in $V$.*

*We say that $X$ is* smoothly locally $\mathfrak{T}$–regular *in $Y$ if every smoothly local base of $X$ in $Y$ is regular.*

The difficulty with this is that when considering inheritance from $X$ to $C^\infty(S, X)$, not every neighbourhood of $C^\infty(S, X)$ in $C^\infty(S, Y)$ is of the form $C^\infty(S, V)$. However, for many applications it is enough to consider those neighbourhoods which are of this form.

**Theorem 3.14** *Let $X$ and $Y$ be smooth spaces with $X \subseteq Y$ and $X$ smoothly locally $\mathfrak{T}$–regular in $Y$. Let $S$ be a smoothly $\mathfrak{T}$–compact space that is also weakly smoothly compact. Then for a smoothly local base $X \overset{o}{\subseteq} V \subseteq Y$, $C^\infty(S, X) \overset{o}{\subseteq} C^\infty(S, V) \subseteq C^\infty(S, Y)$ is regular.*

*Proof.* Let $X \overset{o}{\subseteq} V \subseteq Y$ be a smoothly local base of $X$ in $Y$ and let $f\colon Y \to \mathbb{R}$ be a smooth function which is 1 on $X$ and supported in $V$. Applying the mapping functor $C^\infty(S, -)$ we obtain a smooth map:

$$C^\infty(S, f)\colon C^\infty(S, Y) \to C^\infty(S, \mathbb{R}).$$

This has the following properties:

1. If $\alpha \in C^\infty(S, Y)$ has image in $X$ then $f(\alpha)$ is the constant map $s \mapsto 1$.

2. If $\alpha \in C^\infty(S, Y)$ is such that for all $s \in S$ then $f(\alpha)(s) \ne 0$ then $\alpha \in C^\infty(S, V)$.

As $S$ is weakly smoothly compact there is a smooth function $g\colon C^\infty(S, \mathbb{R}) \to [0, 1]$ such that $g(0) = 1$ and if $g(\beta) \ne 0$ then $\beta \in C^\infty(S, (-1, 1))$. Let $h\colon C^\infty(S, Y) \to \mathbb{R}$ be the composition $\alpha \mapsto g(1 - f(\alpha))$. Then $h$ has the following properties:

1. If $\alpha \in C^\infty(S, Y)$ has image in $X$ then $h(\alpha) = g(1 - 1) = g(0) = 1$.

2. If $\alpha \in C^\infty(S, Y)$ is such that $h(\alpha) \ne 0$ then for $s \in S$, $|1 - f(\alpha)(s)| < 1$ whence $f(\alpha)(s) \ne 0$ and so $\alpha \in C^\infty(S, V)$.

Hence $h$ is a smooth function $C^\infty(S, Y) \to \mathbb{R}$ which is 1 on $C^\infty(S, X)$ and is supported in $C^\infty(S, V)$.

## 3.4 Miscellaneous Properties

There is one more thing that we would like to introduce which fits in with the general theme of this article but which is a little disjoint from the other things under consideration. That is the notion of a fibre bundle.

A smooth map $\pi\colon E \to B$ is a fibre bundle if for each $b \in B$ there is a neighbourhood $V_b$ and a trivialisation $E|_{V_b} \cong V_b \times F$, where $F$ is the model for the fibre. Where we do not have a "standard model" for the fibre, we can use $E_b$ instead. Thus we have a bundle isomorphism $E|_{V_b} \cong V_b \times E_b$. Now we let $b$ vary over $B$. Taking disjoint unions (thus ignoring the smooth structure for the moment), our base space is $\coprod_b V_b$, which we can write as $V := \{(b, p) \in B \times B : p \in V_b\}$. On the left, we have $\coprod_b E|_{V_b}$. The fibre of this at $(b, p)$ is $E_p$, whence this is $\pi_2^* E$. On the right, we have $\coprod_b V_b \times E_b$. The fibre of this at $(b, p)$ is $E_b$, whence this is $\pi_1^* E$. Hence our isomorphisms fit together to a bijection (as we're ignoring smooth structures) $\pi_2^* E \cong \pi_1^* E$. This motivates our definition.

**Definition 3.15** *A smooth map $\pi\colon E \to B$ of smooth spaces is a* smoothly local fibre bundle *if there is an open neighbourhood $V \subseteq B \times B$ of the diagonal and a bundle isomorphism $\pi_1^* E \cong \pi_2^* E$ over $V$.*

We can obviously augment this by putting restrictions on what a "bundle isomorphism" is for special types of bundle: for example, vector bundle and principal bundle. It is also straightforward to say what a trivial bundle looks like in this picture: it is one where we may take $V = B \times B$.

As with all of our properties to date, this is obviously inherited by mapping spaces.

**Proposition 3.16** *Let $\pi\colon E \to B$ be a smoothly local fibre bundle. Let S be a smoothly $\mathfrak{T}$–compact smooth space. Then $C^\infty(S,\pi)\colon C^\infty(S,E) \to C^\infty(S,B)$ is a smoothly local fibre bundle.*

## 4 Tangent Spaces

The standard definition of a local addition for a finite dimensional smooth manifold $M$ gives a diffeomorphism $\varphi\colon TM \to V \subseteq M \times M$ of spaces over $M$. One use of this is to define a smoothly locally additive structure on $M$ by transferring the vector bundle structure from $TM$ to $V$. But this doesn't make any special use of the fact that we start with the tangent bundle.

To see what extra we gain from using the tangent bundle, let us start at the other end. Using the scalar multiplication from a smoothly locally additive structure on $M$ we can define a map $V \to TM$ by sending $v \in V$ to the equivalence class of the curve $t \mapsto \pi_2(tv)$. This is a smooth map of spaces over $M$. Under reasonable conditions on the local addition this is the inverse map. This gives us a way to smoothly choose representing curves for tangent vectors and this is an extremely useful thing to be able to do whence it is the key feature that we would like to replicate.

### 4.1 Defining Tangent Spaces

First, we need to define the tangent space of a smooth space. Any cocomplete cartesian closed category of generalised smooth spaces admits a notion of the tangent space of a smooth space. The ways of defining the tangent space of a smooth manifold that do not depend on a choice of coordinates generalise, but are not necessarily still the same. In this paper, we shall concentrate on the smallest sensible definition: equivalence classes of smooth curves. To distinguish it from the other definitions we shall refer to it, and objects derived from it, as *kinematic*.

**Definition 4.1** *Let $\mathcal{S}$ be a cocomplete cartesian closed category of generalised smooth spaces. Let X be a smooth space.*

*Let $W_X \subseteq C^\infty(\mathbb{R},X) \times C^\infty(\mathbb{R},X)$ be the subspace consisting of pairs $(\alpha,\beta)$ of smooth curves in X satisfying $\alpha(0) = \beta(0)$ and for each $f \in C^\infty(X,\mathbb{R})$ then $(f \circ \alpha)'(0) = (f \circ \beta)'(0)$. There are two obvious maps $W_X \to C^\infty(\mathbb{R},X)$. Define the* kinematic tangent space *of X, written TX, as the coequaliser of $W_X \rightrightarrows C^\infty(\mathbb{R},X)$.*

*For a smooth map $g\colon X \to Y$ we note that post-composition by g defines a smooth map $W_X \to W_Y$ and hence induces a smooth map $Tg\colon TX \to TY$.*

*The* anchor map *is the projection $\pi\colon TX \to X$ induced by evaluation of curves at 0.*

*For a point, say $x \in X$, we define the* tangent space of X at x *to be the fibre of $\pi$ at x considered as a subspace of TX. We write this as $T_xX$.*

As $\mathcal{S}$ is a topological concrete category coequalisers in $\mathcal{S}$ are formed by taking the quotient in Set and then imposing the initial smooth structure for the quotient map from the path space. Thus an element (*tangent vector*) of $TX$ is represented by a smooth path $\alpha\colon \mathbb{R} \to X$ and two paths, say $\alpha$ and $\beta$, represent the same tangent vector if and only if they satisfy $\alpha(0) = \beta(0)$ and for all $f \in C^\infty(X,\mathbb{R})$ then $(f \circ \alpha)'(0) = (f \circ \beta)'(0)$. If smooth functions separate points in $X$ then this relation can be expressed as $(f \circ \alpha)^{(j)}(0) = (f \circ \beta)^{(j)}(0)$ for $j = 0,1$.

It is clear that the assignment $X \to TX$ is a functor.

There is an alternative choice for the smooth structure on the tangent space at a point. We could define it as the quotient of $C^\infty((\mathbb{R}, 0), (X, x))$ by the same equivalence relation. Whether or not these are equivalent structures in general will depend somewhat on how quotients are constructed in the category of smooth objects.

## 4.2 Basic Properties

As we are currently working in an extremely general situation we cannot say definitively that the kinematic tangent space of an ordinary manifold agrees with its standard tangent space. It is certainly the case that the two will have the same (or naturally isomorphic) underlying set and the identity (or the natural equivalence) defines a smooth map from the standard to the kinematic (we will describe this in a moment). What we cannot completely rule out is that the expansion of the category from that of smooth manifolds to that of smooth spaces has introduced some strange smooth spaces in between these.

The map from the standard to the kinematic can be described as follows. Let $M$ be a smooth manifold and write $T^{(s)}M$ for its standard tangent space. Pick a complete Riemannian metric on $M$. Let $\exp \colon \mathbb{R} \times T^{(s)}M \to M$ be the exponential map. As we are in a cartesian closed category this adjuncts to the geodesic map $T^{(s)}M \to C^\infty(\mathbb{R}, M)$ which we can compose with the quotient mapping $C^\infty(\mathbb{R}, M) \to TM$.

Obviously, the best situation is when this natural transformation is a natural isomorphism. As the category of smooth spaces is considered as an *extension* of that of smooth manifolds, we shall use the following terminology as a shortcut for this property.

**Definition 4.2** *Let $\mathcal{S}$ be a suitable category of generalised smooth spaces. We say that $\mathcal{S}$ is a* kinematic extension *of the category of smooth manifolds if the natural transformation from the standard tangent functor on smooth manifolds to the restriction of the kinematic tangent functor is a natural isomorphism.*

As our category is cartesian closed we can ask for whether or not the kinematic tangent functor is *strong* in that the induced map on the hom-sets lifts to a morphism on the internal hom-objects.

**Lemma 4.3** *The kinematic tangent functor is strong.*

*Proof.* We want to show that the set map $C^\infty(X, Y) \to C^\infty(TX, TY)$ induced by the functor lifts to a smooth map. As our category is cartesian closed it is sufficient to show that the following map is smooth:

$$TX \times C^\infty(X, Y) \to TY. \tag{1}$$

Now the smooth map $C^\infty(\mathbb{R}, X) \to TX$ is a coequaliser. As our category is cartesian closed coequalisers commute with products and so the following is also a coequaliser:

$$C^\infty(\mathbb{R}, X) \times C^\infty(X, Y) \to TX \times C^\infty(X, Y).$$

Therefore the map of (1) is smooth if (and only if) the following induced map is smooth:

$$C^\infty(\mathbb{R}, X) \times C^\infty(X, Y) \to TY. \tag{2}$$

On the set level this is the map $(\alpha, g) \mapsto [g \circ \alpha]$ and therefore it factors through $C^\infty(\mathbb{R}, Y) \to TY$. Hence to show that the map in (2) is smooth it is sufficient to show that composition is smooth as a map:

$$C^\infty(\mathbb{R}, X) \times C^\infty(X, Y) \to C^\infty(\mathbb{R}, Y)$$

which is a standard property of cartesian closed categories.

Hence the map $C^\infty(X, Y) \to C^\infty(TX, TY)$ induced by the kinematic tangent functor is smooth and thus the functor is strong.

Now let us consider the interaction between tangent spaces and products. As tangent spaces are defined as quotients this interaction is not guaranteed to be nice. However, under reasonable assumptions on the category we can show good behaviour.

**Proposition 4.4** *Let $\Lambda$ be an indexing set. Assume that $\mathbb{R}^\Lambda$ has the property that for $f \colon \mathbb{R}^\Lambda \to \mathbb{R}$ with $Tf$ non-zero at $0$ then there is some direction $\lambda \in \Lambda$ which detects this.*

*Let $\{X_\lambda\}_{\lambda \in \Lambda}$ be a family of smooth spaces indexed by $\Lambda$. There is a natural map $T(\prod_\lambda X_\lambda) \to \prod_\lambda TX_\lambda$ coming from the evaluation maps $\prod_\lambda X_\lambda \to X_\lambda$ which is a surjection on the underlying sets. Under the above assumption it is also injective (whence a bijection).*

*In particular, this is true whenever $\Lambda$ is finite.*

*Proof.* To see that it is surjective we take a family of tangent vectors $v_\lambda \in TX_\lambda$. For each, we choose a representing smooth curve $\alpha_\lambda \colon \mathbb{R} \to X_\lambda$. Then the curve $\alpha := (\alpha_\lambda) \colon \mathbb{R} \to \prod_\lambda X_\lambda$ is smooth and so represents a tangent vector in $T \prod_\lambda X_\lambda$. By construction, the equivalence class of $\alpha$ maps to the family of equivalence classes of the $\alpha_\lambda$ and hence $T \prod_\lambda X_\lambda \to \prod_\lambda TX_\lambda$ is surjective.

Now we assume the condition of the proposition. Let $\alpha, \beta \colon \mathbb{R} \to \prod_\lambda X_\lambda$ be smooth curves with $\alpha(0) = \beta(0)$ and assume that there is some $f \colon \prod_\lambda X_\lambda \to \mathbb{R}$ with the property that $(f \circ \alpha)'(0) \neq (f \circ \beta)'(0)$. For $\lambda \in \Lambda$, let $\alpha_\lambda \colon \mathbb{R} \to X_\lambda$ be the composition of $\alpha$ with the projection to $X_\lambda$, and $\beta_\lambda$ likewise. These provide a smooth map $\prod_\lambda \alpha_\lambda \colon \prod_\lambda \mathbb{R} \to \prod_\lambda X_\lambda$ (and similarly for $\beta$) with the property that after setting $\Delta \colon \mathbb{R} \to \prod_\lambda \mathbb{R}$ to be the diagonal we have $\alpha = \prod_\lambda \alpha_\lambda \circ \Delta$. Now consider the two maps $f_\alpha := f \circ \prod_\lambda \alpha_\lambda$ and $f_\beta := f \circ \prod_\lambda \beta_\lambda$. As we are mapping into $\mathbb{R}$ we can consider their difference, $f_\alpha - f_\beta$. By construction $T(f_\alpha \circ \Delta)$ and $T(f_\beta \circ \Delta)$ differ at $0$. Thus $T(f_\alpha - f_\beta)$ must be non-zero at $\prod_\lambda 0$. The assumption now says that there is some direction, say $\lambda_0$, which detects this. Let $\iota_{\lambda_0} \colon \mathbb{R} \to \mathbb{R}^\Lambda$ be the linear inclusion of the $\lambda_0$ direction. Then the conclusion is that $T(f_\alpha \circ \iota_{\lambda_0}) \neq T(f_\beta \circ \iota_{\lambda_0})$. Tracing through the construction, this implies that $\alpha_{\lambda_0}$ and $\beta_{\lambda_0}$ represent different tangent vectors in $X_{\lambda_0}$. Hence $T \prod_\lambda X_\lambda \to \prod_\lambda TX_\lambda$ is injective.

Let us consider the situation when $\Lambda$ is finite, thus $\mathbb{R}^\Lambda$ is $\mathbb{R}^n$ for some finite $n$. Then we have a smooth map $\mathbb{R}^{2n} \to T\mathbb{R}^n$ by sending $(u, v)$ to the equivalence class of the curve $t \mapsto u + tv$. This is a bijection on the underlying sets as quotients in $\mathcal{S}$ are formed by taking the quotient of sets and then putting on the quotient smooth structure, and the construction of $T\mathbb{R}^n$ from $C^\infty(\mathbb{R}, \mathbb{R}^n)$ and $C^\infty(\mathbb{R}^n, \mathbb{R})$ gives the ordinary tangent space of $\mathbb{R}^n$ as a set. This is the natural transformation on the full subcategory of Euclidean spaces from the "standard" to the "kinematic" tangent space functors.

Let $f \colon \mathbb{R}^n \to \mathbb{R}$ be a smooth map. The above natural transformation intertwines the tangent map $Tf \colon T\mathbb{R}^n \to T\mathbb{R}$ with the differential map $(f, Df) \colon \mathbb{R}^{2n} \to$ $\mathbb{R}^2$. So if $Tf$ is non-zero at $0 \in \mathbb{R}^n$, $Df$ is likewise non-zero and thus there is a direction, say $j$, which detects this. That is to say, if $x^j \colon \mathbb{R} \to \mathbb{R}^n$ is the insertion of the $j$th coordinate then $f \circ x$ has non-trivial derivative at $0$. This is what is required for the statement of the Proposition to hold. 🐦

**Corollary 4.5** *The kinematic tangent functor commutes with finite products.*

*Proof.* As our category is cartesian closed taking the product with a space preserves coequalisers. Thus for smooth spaces $X$ and $Y$ the maps $C^\infty(\mathbb{R}, X) \times C^\infty(\mathbb{R}, Y) \to C^\infty(\mathbb{R}, X) \times TY$ and $C^\infty(\mathbb{R}, X) \times TY \to TX \times TY$ are coequalisers. Since we're in a topological category the composition of coequalisers is again a coequaliser. Hence the map $C^\infty(\mathbb{R}, X) \times C^\infty(\mathbb{R}, Y) \to TX \times TY$ is a coequaliser. The natural map $T(X \times Y) \to TX \times TY$ is thus a map of coequalisers. From Proposition 4.4 it is a bijection on underlying sets and thus an isomorphism of smooth spaces. 🐦

Another situation where tangent spaces might not behave quite as desired is with subspaces. In particular, we are interested in the situation of a $\mathfrak{T}$–open subspace of an ambient space.

**Lemma 4.6** *Let $X$ be a smooth space. Let $x \in X$ and let $U \subseteq X$ be a neighbourhood of $x$. Then $T_x U \to T_x X$ is surjective.*

*Suppose that there is a neighbourhood $W$ of $x$ with $W \subseteq V$ and the property that the images of the restriction maps $C^\infty(X, \mathbb{R}) \to C^\infty(W, \mathbb{R})$ and $C^\infty(U, \mathbb{R}) \to C^\infty(W, \mathbb{R})$ coincide. Then $T_x U \to T_x X$ is injective.*

*Proof.* Let $v \in T_x X$. Then $v$ is represented by a smooth curve $c \colon \mathbb{R} \to X$ with $c(0) = x$. As $U$ is a neighbourhood of $x$ there is some $\varepsilon > 0$ such that for $|t| < \varepsilon$ then $c(t) \in U$. Let $a \colon \mathbb{R} \to \mathbb{R}$ be a smooth function with $a(0) = 0$, $a'(0) = 1$, and for all $t \in \mathbb{R}$, $|a(t)| < \varepsilon$. Consider the curve $c \circ a \colon \mathbb{R} \to X$. Let $f \colon X \to \mathbb{R}$ be a smooth function. Then $f \circ (c \circ a) = (f \circ c) \circ a$. Since $f \circ c$ and $a$ are smooth functions $\mathbb{R} \to \mathbb{R}$, the chain rule implies that $(f \circ (c \circ a))'(0) = (f \circ c)'(0)a'(0)$ (using $a(0) = 0$), whence $(f \circ (c \circ a))'(0) = (f \circ c)'(0)$. Hence $c \circ a$ also represents $v$. Now $c \circ a$ factors through $U$ and hence represents a tangent vector in $T_x U$. Thus $T_x U \to T_x X$ is surjective.

The reason that it can fail to be injective is because there may be more functions $U \to \mathbb{R}$ than just those that are restrictions of functions on $X$. The existence of the smaller neighbourhood ensures that these don't change the tangent space at $x$.

Let $v_1, v_2 \in T_x U$ be distinct tangent vectors at $x$. As in the above, they can be represented by curves $c_1, c_2 \colon \mathbb{R} \to U$ with the property that their images lie wholly in $W$. As they are distinct there is some $f \colon U \to \mathbb{R}$ that distinguishes them. That is, that $(f \circ c_1)'(0) \neq (f \circ c_2)'(0)$. By assumption, there is some $g \colon X \to \mathbb{R}$ such that $g|_W = f|_W$. Then $(g \circ c_1)'(0) \neq (g \circ c_2)'(0)$ and thus $c_1$ and $c_2$ represent different tangent vectors in $T_x X$, whence $T_x U \to T_x X$ is injective. 🐦

The condition for injectivity occurs, for example, when the smooth functionals on a smooth space obey a sheaf-like condition. More precisely, we have the following result.

**Proposition 4.7** *Let $X$ be a smooth space. For $f \in C^\infty(X, \mathbb{R})$ let $U_f := \{x \in X : f(x) \neq 0\}$ and $V_f := \{x \in X : f(x) \neq 1\}$.*

*We assume that for every smooth $f \colon X \to \mathbb{R}$, there is a pull-back diagram:*

$$\begin{array}{ccc} C^\infty(X,\mathbb{R}) & \longrightarrow & C^\infty(U_f,\mathbb{R}) \\ \downarrow & \lrcorner & \downarrow \\ C^\infty(V_f,\mathbb{R}) & \longrightarrow & C^\infty(U_f \cap V_f,\mathbb{R}) \end{array}$$

*Then for $f \in C^\infty(X,\mathbb{R})$ and $x \in U_f$, $T_x U_f \to T_x X$ is a bijection.*

*Proof.* Let us start by replacing $f$ by a more convenient function, which we shall also call $f$. We divide by $f(x)$ to ensure that it has value 1 at $x$, then square to ensure that the function is positive. None of this changes $U_f$.

Let $h_1 \colon X \to \mathbb{R}$ be the composition of $f$ with a smooth function $\mathbb{R} \to \mathbb{R}$ with image in $[0,1]$ which maps $[{}^3\!/\!_4, \infty)$ to $\{1\}$ and $(-\infty, {}^1\!/\!_4)$ to $\{0\}$. Similarly, let $h_2 \colon X \to \mathbb{R}$ be the composition of $f$ with a smooth function $\mathbb{R} \to \mathbb{R}$ with image in $[0,1]$ which maps 1 to 1 and $(-\infty, {}^3\!/\!_4)$ to $\{0\}$. Note that $U_{h_1} \subseteq U_f$, that $h_2(x) = 1$, and that if $h_1(y) \neq 1$ then $h_2(y) = 0$.

Let $g \in C^\infty(U_f, \mathbb{R})$ and let $g_0 \in C^\infty(U_{h_1}, \mathbb{R})$ be the restriction of $h_2 g$ to $U_{h_1}$, which is defined since $U_{h_1} \subseteq U_f$. Let $g_1 \in C^\infty(V_{h_1}, \mathbb{R})$ be the zero function. Now for $y \in U_{h_1} \cap V_{h_1}$, $h_1(y) \neq 1$ whence $h_2(y) = 0$. Hence $(h_2 g)(y) = 0$ and so $g_0$ restricts to the zero function on $U_{h_1} \cap V_{h_1}$. Thus, by assumption, there is a function $\hat{g} \colon X \to \mathbb{R}$ which restricts to $g_0$ on $U_{h_1}$. It is clear that this also restricts to $h_2 g$ on $U_f$.

Thus the conditions for the second part of Lemma 4.6 are full-filled and $T_x U_f \to T_x X$ is injective. ☙

This condition holds, for example, for all of the categories of generalised smooth spaces considered in [Sta11].

Although we refer to elements of the tangent space as *tangent vectors* the fibres might not be vector spaces. The following proposition shows what structure they do have and how far they are from being vector spaces.

**Proposition 4.8** *Let $X$ be a smooth space. Let us say that for $x \in X$ and tangent vectors $u, v \in T_x X$ then $u$ and $v$ are* summable *if there is some $w \in T_x X$, which we shall refer to as their sum, such that for all smooth functions $f \colon X \to \mathbb{R}$ then $Tf(w) = Tf(u) + Tf(v)$.*

*The tangent space, $TX$, has the following structure:*

1. *A section $0 \colon X \to TX$ defined by $0(x) = [t \mapsto x]$.*

2. *An action of $\mathbb{R}$, $\lambda \colon \mathbb{R} \times TX \to TX$, defined by $\lambda(c, [\alpha]) = [t \mapsto \alpha(ct)]$.*

3. *If $u, v$ are summable, their sum is unique.*

4. *The zero, $\mathbb{R}$–action, and summation obey all of the identities of a vector space over $X$ whenever all the terms are defined.*

*Proof.* The claims for the zero section and scalar multiplication are obvious.

To show that the sum is unique, let $\alpha$ and $\beta$ represent $u$ and $v$. Let $w_1$ and $w_2$ be two sums of $u$ and $v$ and let $\gamma_1$ and $\gamma_2$ represent them. By definition, we have $\gamma_1(0) = \gamma_2(0) = \alpha(0) (= \beta(0))$. Let $f \colon X \to \mathbb{R}$ be a smooth map. Then under the identification of $T\mathbb{R}$ with $\mathbb{R} \times \mathbb{R}$, we have $(f \circ \gamma_1, (f \circ \gamma_1)')(0) = Tf(w_1)$. By assumption, $Tf(w_1) = Tf(u) + Tf(v) = Tf(w_2)$ and thus $(f \circ \gamma_1)'(0) = (f \circ \gamma_2)'(0)$. Hence $\gamma_1$ and $\gamma_2$ represent the same tangent vector in $TX$, whence $w_1 = w_2$.

The identities for a vector space structure (over X) are proved in a similar fashion. 🕊

In synthetic differential geometry smooth spaces where the tangent fibres are vector spaces are called *microlinear*.

## 4.3 Smoothly Locally Kinematic Spaces

Now that we have the notion of tangent spaces we can define a smoothly locally kinematic space. It builds on the definition of a smoothly locally deflatable space. Let $X$ and $Y$ be smooth spaces and $X \mathrel{\widehat{\subseteq}} V \subseteq Y$ be a smoothly local base for $X$ in $Y$, with retraction $\tau$. Let $H \colon \mathbb{R} \times V \to V$ be a smoothly local deflation. Define $\ker T\tau \subseteq TV$ to be the subspace of tangent vectors which map to the corresponding zero vector under the mapping $T\tau$. Let $\iota \colon X \to V$ be the inclusion and define $\iota^* \ker T\tau$ to be pull-back in the following diagram.

$$\begin{array}{ccc} \iota^* \ker T\tau & \longrightarrow & \ker T\tau \\ \downarrow & & \downarrow \pi \\ X & \stackrel{\iota}{\longrightarrow} & V \end{array}$$

Define a smooth map $\Psi \colon V \to C^\infty(\mathbb{R}, V)$ by $v \mapsto (t \mapsto H_t(v))$. Then we have $\Psi(v)(0) = \tau(0)$ so the resulting curves are anchored at points of $X$. In addition, $\tau H_t(v) = \tau(v)$ so the tangent vector represented by $t \mapsto H_t(v)$ lies in $\ker T\tau$. These conditions combine to show that $\Psi$ defines a smooth map $\psi \colon V \to \iota^* \ker T\tau$.

**Definition 4.9** *Let $X \subseteq Y$ be smooth spaces with $X$ smoothly locally deflatable in $Y$. Let $X \mathrel{\widehat{\subseteq}} V \subseteq Y$ be a smoothly local base admitting a smoothly local deflation $H \colon \mathbb{R} \times V \to V$. We define the* kinematic map *associated to this deflation to be the smooth map $\psi \colon V \to \iota^* \ker T\tau$ given by the above construction.*

**Definition 4.10** *Let $X \subseteq Y$ be smooth spaces. A smoothly locally deflatable structure on $X$ in $Y$ is* kinematic *if the kinematic map $\psi \colon V \to \iota^* \ker T\tau$ is a diffeomorphism.*

*We say that $X$ is* smoothly locally kinematic in $Y$ *if it admits a kinematic smoothly locally deflatable structure.*

### 4.3.1 Kinematic Vector Spaces

The simplest example to start with is that of a smooth vector space, say $E$. We can define a smoothly locally deflatable structure on $E$ by scalar multiplication. Explicitly, let $V = E \times E$ with $\tau \colon V \to E$ projection onto the first coordinate. Define $H \colon \mathbb{R} \times E \times E \to E \times E$ by $H_t(e_1, e_2) = (e_1, e_1 + t(e_2 - e_1))$. By Corollary 4.5, $T(E \times E) \cong TE \times TE$, whence $\iota^* \ker T\tau \cong TE$. Then if the kinematic map $\psi \colon V \to \iota^* \ker T\tau$ is a diffeomorphism every tangent vector in $TE$ can be represented by a curve of the form $t \mapsto e_1 + te_2$ for some $e_2 \in E$. That is to say, the natural map $E \times E \to TE$ is a diffeomorphism. Since $TE$ splits as $E \times T_0 E$, we can simplify this to the natural linear map $E \to T_0 E$ and make the following definition.

**Definition 4.11** *Let $E$ be a smooth vector space. We say that $E$ is a* kinematic *vector space if the canonical map $E \to T_0 E$ is a diffeomorphism.*

### 4.3.2 Kinematic and Additive

Let us consider the case where the smooth space $X$ is smoothly locally kinematic. Interpreting the definition we see that this means that there is a neighbourhood, say $V$, of the diagonal in $X \times X$ with a smoothly locally deflatable structure such that the kinematic map $\psi \colon V \to \iota^* \ker T\tau$ is a diffeomorphism.

Further structure is provided by the next result which provides the basis for showing that under certain circumstances then a smoothly locally kinematic space is smoothly locally additive.

**Proposition 4.12** *Let $X$ be a smoothly locally deflatable space with structure $X \overset{\frown}{\subseteq} V \subseteq X \times X$, $\tau \colon V \to X$, and let $\psi \colon V \to \iota^* \ker T\tau$ be the kinematic map.*

*Then there is an open set $W \subseteq \mathbb{R}^2 \times V \times_X V$ and a smooth map $\sigma \colon W \to V$ with the following properties:*

1. *$\{(0,0)\} \times V \times_X V \subseteq W$,*

2. *the restriction of $\sigma$ to $\{(0,0)\} \times V \times_X V$ is the projection to $X$ (viewed as a subspace of $V$),*

3. *let $\partial_x$ at $(0,0,u,v) \in W$ be the tangent vector in the first direction and $\partial_y$ in the second, then $T\sigma(\partial_y) = \psi(u)$ and $T\sigma(\partial_x) = \psi(v)$.*

*Proof.* Normally we are a bit loose about identifying $X$ with the diagonal in $X \times X$. Here we want to be more careful. So the retraction $\tau$ is actually the function $V \to \Delta(X)$ given by $\tau = \Delta \circ \pi_1$.

Let $V \times_X V$ be the space of pairs of points in $V$ that lie over the same point in $X$. Define a smooth map $\gamma \colon \mathbb{R} \times V \times_X V \to X \times X$ by

$$\gamma(t,u,v) = (\pi_2 H_t(u), \pi_2(v)).$$

Let $W_\gamma = \gamma^{-1}(V)$. At $t = 0$ we have $H_0(u) = \tau(u) = \Delta \circ \pi_1(u)$ so $\pi_2 H_0(u) = \pi_1(u) = \pi_1(v)$ and thus $\gamma(0,u,v) = v$. Hence $\{0\} \times V \times_X V \subseteq W_\gamma$.

Now define $\sigma \colon \mathbb{R} \times W_\gamma \to X \times X$ by:

$$\sigma(s,t,u,v) = (\pi_1 v, \pi_2 H_s \gamma(t,u,v)).$$

Let $W_\sigma = \sigma^{-1}(V)$. At $s = t = 0$ we have:

$$\sigma(0,0,u,v) = (\pi_1 v, \pi_2 H_0(0,u,v)) = (\pi_1 v, \pi_2 H_0 v)$$
$$= (\pi_1 v, \pi_2 \tau v) = (\pi_1 v, \pi_1 v) = \tau v$$

whence $\sigma(0,0,u,v)$ is on the diagonal at $\tau v = \tau u$. Hence $\{(0,0)\} \times V \times_X V \subseteq W_\sigma$ and the restriction of $\sigma$ to this subspace is the projection to $\Delta(X)$.

Consider the curve defined by $\sigma$ with $s = 0$ and for $t$ small:

$$\begin{aligned}
\sigma(0,t,u,v) &= (\pi_1 v, \pi_2 H_0 \gamma(t,u,v)) \\
&= (\pi_1 v, \pi_2 \tau \gamma(t,u,v)) & H_0 &= \tau \\
&= (\pi_1 v, \pi_1 \gamma(t,u,v)) & \pi_2 \tau &= \pi_1 \\
&= (\pi_1 v, \pi_2 H_t(u)) \\
&= (\pi_1 H_t(u), \pi_2 H_t(u)) & \pi_1 v &= \pi_1 u = \pi_1 H_t(u) \\
&= H_t(u).
\end{aligned}$$

Thus $t \mapsto \sigma(0, t, u, v)$ is (for small $t$) $H_t(u)$.

Now consider the curve defined by $\sigma$ with $t = 0$ and for $s$ small:

$$\begin{aligned}
\sigma(s, 0, u, v) &= (\pi_1 v, \pi_2 H_s \gamma(0, u, v)) & \\
&= (\pi_1 v, \pi_2 H_s(v)) & \gamma(0, u, v) = v \\
&= (\pi_1 H_s(v), \pi_2 H_s(v)) & \pi_1 H_s(v) = \pi_1(v) \\
&= H_s(v).
\end{aligned}$$

Thus $s \mapsto \sigma(s, 0, u, v)$ is $H_s(v)$.

Hence at $(0, 0, u, v)$, $T\sigma$ maps $\partial_x$ to $\psi(v)$ and $\partial_y$ to $\psi(u)$.

**Corollary 4.13** *Let $X$ be a smoothly locally deflatable space. Then the fibres of $\psi(V) \subseteq \iota^* \ker T\tau$ are vector spaces.*

*If $X$ is smoothly locally kinematic then the fibres of $\iota^* \ker T\tau$ and the fibrewise tangent spaces of $X$ are vector spaces.*

*Proof.* We assume the notation of Proposition 4.12 and its proof.

In both cases the thing to prove is the existence of a sum of two vectors, the other properties will follow from Proposition 4.8; in the case of $\iota^* \ker T\tau$ this is because it embeds in $TV$.

For $x \in X$ then the case of $T_x X$ follows from that of $\iota^* \ker T\tau$ at $x$ by surjectivity of the map $\iota^* \ker T\tau \to TX$.

Let $x \in X$ and let $u, v$ be in the fibre of $V$ at $x$. Then $\sigma \colon W \to V$ has the property that $T\sigma(\partial_x) = \psi(v)$ and $T\sigma(\partial_y) = \psi(u)$. Hence $T\sigma(\partial_x + \partial_y)$ exists in $\iota^* \ker T\tau$ and is the sum of $\psi(u)$ and $\psi(v)$.

The important question is as to the smoothness of this construction for a smoothly locally kinematic space. For $\iota^* \ker T\tau$ the main question is whether the pointwise tangent vectors $\partial_x$ and $\partial_y$ can be chosen smoothly. For $TX$ there is the additional question as to whether the map $\iota^* \ker T\tau \to TX$ has a smooth retraction (or, preferably, is a diffeomorphism).

### 4.3.3 Mapping Spaces

The obvious question to ask of smoothly locally kinematic spaces is whether or not the mapping space of one is again smoothly locally kinematic. Since the mapping space construction preserves smoothly locally deflatable structures, the kinematic map $V \to \iota^* \ker T\tau$ exponentiates to the corresponding map:

$$C^\infty(S, V) \to C^\infty(S, \iota^* \ker T\tau).$$

Now $\iota^* \ker T\tau$ is defined as a particular subspace of $TV$. We want to compare this with $\iota^* \ker TC^\infty(S, \tau)$ (where, by abuse of notation, we have used $\iota$ again for $C^\infty(S, \iota)$).

A useful thing here is the natural transformation from $TC^\infty(S, -)$ to $C^\infty(S, T-)$ which at $s \in S$ is the tangent of the evaluation map. This comes from applying the tangent functor to the counit of the adjunction $S \times C^\infty(S, X) \to X$ together with the isomorphism from Corollary 4.5 and the zero section $S \to TS$ to get $S \times TC^\infty(S, X) \to TX$ whence a smooth map $TC^\infty(S, X) \to C^\infty(S, TX)$.

Now consider a smoothly locally deflatable structure on $X \overset{\frown}{\subseteq} V \subseteq Y$ with retraction $\tau$ and kinematic map $\psi \colon V \to \iota^* \ker T\tau$. This induces a smoothly locally deflatable structure on the mapping spaces as well. We have the following commutative diagram.

$$
\begin{array}{ccccc}
C^\infty(S,V) & \xrightarrow{\psi_{C^\infty(S,-)}} & \iota^* \ker TC^\infty(S,\pi) & \longrightarrow & TC^\infty(S,V) \\
\downarrow \text{Id} & & \downarrow & & \downarrow \\
C^\infty(S,V) & \xrightarrow{C^\infty(S,\psi)} & C^\infty(S, \iota^* \ker T\pi) & \longrightarrow & C^\infty(S,TV)
\end{array}
$$

From this diagram it is clear that the induced smoothly locally deflatable structure on the mapping space is again smoothly locally kinematic if (and only if) the map $\iota^* \ker TC^\infty(S,\pi) \to C^\infty(S, \iota^* \ker T\pi)$ is a diffeomorphism. The smoothly locally kinematic structure on $X \subseteq Y$ defines a splitting of this map via

$$\psi_{C^\infty(S,-)} \circ C^\infty(S,\psi)^{-1}.$$

The problem is, therefore, to determine when this map is injective. If it is, the mapping space is smoothly locally kinematic.

## 4.4 Fermat Spaces

The question at the end of the previous section pertains to the natural map $TC^\infty(S,X) \to C^\infty(S,TX)$. This is related to the situation in Proposition 4.4. In particular, if $S$ is discrete then Proposition 4.4 gives conditions where $TC^\infty(S,X) \to C^\infty(S,TX)$ is injective. That situation does not generalise to arbitrary $S$ as it relies on being able to consider the points of the source space independently. To work with more general $S$ we need to bring in another tool from ordinary calculus.

As we noted in Section 4.2, it need not be the case that $T\mathbb{R}^n \cong \mathbb{R}^{2n}$ but that this is a desirable situation to be in. If this is the case then the fact that the kinematic tangent functor is strong means that the map $f \mapsto f'$ is smooth as a map $C^\infty(\mathbb{R}, \mathbb{R}^n) \to C^\infty(\mathbb{R}, \mathbb{R}^n)$. On the set level this has a partial inverse given by integration: $g \mapsto \int_0^s g(t) dt$. Again, without knowing the exact smooth structure on $C^\infty(\mathbb{R}, \mathbb{R}^n)$ we cannot conclude that this is smooth. The reason that we would like to be able to do so is to improve on Hadamard's lemma. Recall that Hadamard's lemma (which was an extension of a result of Fermat on polynomials) states that for a smooth function $f \colon \mathbb{R} \to \mathbb{R}^n$ there is a smooth function $g \colon \mathbb{R}^2 \to \mathbb{R}^n$ with the property that:

$$f(x+y) = f(x) + yg(x,y). \tag{3}$$

As $g$ is smooth (in the usual sense) its value at $(x,0)$ must be $f'(x)$. The function $g(x,y)$ is defined as:

$$g(x,y) = \int_0^1 f'(x+sy) ds. \tag{4}$$

We would like to have a parametrised version of this. That is, for a smooth space $S$ and a smooth function $f \colon S \times \mathbb{R} \to \mathbb{R}^n$ we would like to define a smooth function $g \colon S \times \mathbb{R}^2 \to \mathbb{R}^n$ such that:

$$f(s, x+y) = f(s,x) + yg(s,x,y).$$

If $S$ is a smooth manifold then, as we have assumed that the category of smooth manifolds embeds fully and faithfully in our category of smooth spaces, the above formula for $g$ generalises. However, if we want this to work for all smooth spaces $S$ we need to assume that integration is smooth $C^\infty(\mathbb{R}, \mathbb{R}^n) \to C^\infty(\mathbb{R}, \mathbb{R}^n)$.

We would like to go even further than this. We can rearrange Equation (3) to the following form:
$$f(y) - f(x) = y g(x, y). \tag{5}$$
The left-hand side of this should be interpreted as "$f(y)$ as seen from $f(x)$". That is, the subtraction is a relocation rather than an inverse addition. With that interpretation Equation 5 makes sense in a smoothly locally deflatable space.

So let $X$ be a smoothly locally deflatable space with smoothly local base $X \stackrel{\frown}{\subseteq} V \subseteq X \times X$ and deflation $H \colon \mathbb{R} \times V \to V$. As before, we write $H_t$ for $H(t, -)$. In this scenario the statement "$f(y)$ as seen from $f(x)$" is reinterpreted to mean "the pair $(f(x), f(y))$".

This means that instead of $g(x, y)$ being defined for all $x$ and $y$ we can only define it for $x$ and $y$ with $(f(x), f(y)) \in V$. More precisely, let $\eta \colon \mathbb{R}^2 \times C^\infty(\mathbb{R}, X) \to X \times X$ be the double evaluation map, $\eta(s, t, \alpha) = (\alpha(s), \alpha(t))$. Let $V_\eta = \eta^{-1}(V)$. Then we want to be able to say that there is a smooth function $F \colon V_\eta \to V$ with the property that for $(s, t, \alpha) \in V_\eta$:
$$(\alpha(s), \alpha(t)) = H_{t-s} F(s, t, \alpha). \tag{6}$$

For $s \ne t$ then Equation 6 determines $F(s, t, \alpha)$. In the case of $\mathbb{R}$, then $g(x, x)$ was $f'(x)$. In the general case $\alpha'(s)$ is a vector in the tangent space of $X$ at $\alpha(s)$. Thus to compare $F(s, s, \alpha)$ with $\alpha'(s)$ we need to transport it from $V$ to $TX$. The obvious way to do this is via the kinematic map $\psi$ derived from the locally deflatable structure. This leads us to make the following definition.

**Definition 4.14** *Let $X$ be a smooth space. Let $X \stackrel{\frown}{\subseteq} V \subseteq X \times X$ be a smoothly local base and suppose that it admits a smoothly local deflation $H \colon \mathbb{R} \times V \to V$. Let $\psi \colon V \to \iota^* \ker T\tau$ be the kinematic map. Let $\eta \colon \mathbb{R}^2 \times C^\infty(\mathbb{R}, X) \to X \times X$ be the double evaluation map and let $V_\eta = \eta^{-1}(V)$.*

*We say that the smoothly locally deflatable structure on $X$ is* Fermat *if it admits a smooth function $F \colon V_\eta \to V$, which we shall call a* Fermat map, *such that:*

1. *For $s \ne t$ then $\eta(s, t, \alpha) = H_{t-s} F(s, t, \alpha)$,*

2. *$T\pi_2 \circ \psi F(s, s, \alpha) = \alpha'(s)$ (where we think of $T\pi_2$ as a map $\iota^* \ker T\tau \to TX$ via $\iota^* \ker T\tau \hookrightarrow TV \to TX$).*

*We say that $X$ is a* Fermat space *if it admits a smoothly locally deflatable structure that is Fermat.*

Note that if a smoothly locally deflatable structure on $X$ is Fermat then the composition $V \to \iota^* \ker T\tau \to TX$ is surjective.

Any two Fermat maps must agree away from the diagonal, whence if $V$ is functionally Hausdorff (smooth functionals separate points) then there can be at most one Fermat map. This is because on a functionally Hausdorff space a smooth curve is determined by its values on a dense subset.

The main reason for introducing these spaces is the following result.

**Theorem 4.15** *Let X be a Fermat space. Assume that for a suitable choice of deflation, the map $T\pi_2 \circ \psi \colon V \to TX$ is injective where $\psi$ is the associated kinematic map. Let S be a smoothly $\mathfrak{T}$–compact space. Then the kernel of the natural map $TC^\infty(S, X) \to C^\infty(S, TX)$ is null. Equivalently, if $v \in TC^\infty(S, X)$ is non-zero then there is some direction in S which detects that.*

*Proof.* We assume the notation of Definition 4.14.

Let $\alpha \colon \mathbb{R} \to C^\infty(S, X)$ be a smooth curve. Let us write $\hat{\alpha} \colon S \to C^\infty(\mathbb{R}, X)$ for the smooth map defined by $\alpha$. Let $q \colon C^\infty(\mathbb{R}, X) \to TX$ be the quotient and $\zeta \colon X \to TX$ the zero section.

We assume that as smooth maps $S \to TX$ then $q\hat{\alpha} = \zeta\alpha$. That is to say, for $s \in S$ then $\hat{\alpha}(s)'(0) = 0$. Let $f \colon C^\infty(S, X) \to \mathbb{R}$ be a smooth function. We want to show that $(f \circ \alpha)'(0) = 0$.

As S is smoothly $\mathfrak{T}$–compact there is some open interval containing 0 with the property that for all $s \in S$ and for $t$ in that interval, $(\hat{\alpha}(s)(0), \hat{\alpha}(s)(t)) \in V$. We reparametrise $\alpha$ so that we can take the interval to be the whole of $\mathbb{R}$ but without changing $\alpha$ near 0. This reparametrised curve has the same value for $(f \circ \alpha)'(0)$.

Thus for all $t \in \mathbb{R}$ and $s \in S$ then $(\hat{\alpha}(s)(0), \hat{\alpha}(s)(t)) \in V$. So $(s, t) \mapsto (0, t, \hat{\alpha}(s))$ is a smooth map $S \times \mathbb{R} \to V_\eta$ (using the notation of Definition 4.14). We compose this with $F \colon V_\eta \to V$ to get a smooth map $\beta \colon S \times \mathbb{R} \to V$.

For $s \in S$ and $t \in \mathbb{R}$ then by construction $\beta(s, t) = F(0, t, \hat{\alpha}(s))$. Thus $\pi_1\beta(s, t) = \hat{\alpha}(s)(0)$ and for $t \neq 0$ then $\hat{\alpha}(s)(t) = \pi_2 H_t \beta(s, t)$. In addition, $\beta(s, 0)$ is such that $T\pi_2\psi(\beta(s, 0)) = \hat{\alpha}(s)'(0)$. By assumption this is zero and so as $T\pi_2\psi$ is assumed to be injective we must have $\beta(s, 0) = (\hat{\alpha}(s)(0), \hat{\alpha}(s)(0))$.

Now let $\sigma \colon \mathbb{R}^2 \to C^\infty(S, X)$ be the function $\sigma(r, t)(s) = \pi_2 H_r \beta(s, t)$, so that $\alpha(t) = \sigma(t, t)$. We consider the function $f \circ \sigma \colon \mathbb{R}^2 \to \mathbb{R}$. To show that $(f \circ \alpha)'(0) = 0$ we will show that $D(f \circ \sigma)$ vanishes at the origin.

Let us start with $\partial_r(f \circ \sigma)$. This is the derivative of the map $r \mapsto (f \circ \sigma)(r, 0)$. We have $\sigma(r, 0)(s) = \pi_2 H_r \beta(s, 0)$. Then $\beta(s, 0) = (\hat{\alpha}(s)(0), \hat{\alpha}(s)(0))$ which is invariant under $H_r$. Thus $\sigma(r, 0)(s) = \hat{\alpha}(s)(0)$, whence $\partial_r(f \circ \sigma) = 0$.

Let us consider $\partial_t(f \circ \sigma)$. This is the derivative of the map $t \mapsto (f \circ \sigma)(0, t)$. We have $\sigma(0, t)(s) = \pi_2 H_0 \beta(s, t) = \pi_1 \beta(s, t) = \hat{\alpha}(s)(0)$. Again, this is constant.

Hence both partial derivatives of $(f \circ \sigma)$ vanish and so in particular $(f \circ \alpha)'(0) = 0$. As this is true for all $f$ we have $\alpha'(0) = 0$ in $TC^\infty(S, X)$.

## 5 Manifolds

In this section we wish to take the results of the previous sections and apply them to manifolds. In essence, we want to show that ordinary smooth manifolds are about as good as can be. Most of the results in this section are variants on the tubular neighbourhood theorem. Indeed, the weakest results are a simple corollary of the existence of tubular neighbourhoods and need no more technical work. However the more complicated results need tubular neighbourhoods with particular properties and so for those we need to construct them "from scratch" (or, more explicitly, from the inverse function theorem).

We shall assume throughout this section that the category of generalised smooth spaces is a kinematic extension of the category of smooth manifolds, as in Definition 4.2, whence the kinematic tangent space functor restricts to the ordinary tangent space functor for a smooth finite dimensional manifold. In technical terms this is more about notation than anything else. Even if this were not true, the ordinary tangent space of a smooth manifold would still exist as a smooth object and we would use that in place of the kinematic tangent space in the following. All we would lose was the ability to claim that certain spaces were locally kinematic and have to settle for mere locally additive.

**Theorem 5.1** *Let $M$ and $N$ be finite dimensional smooth manifolds. Let $f \colon M \to N$ be a regular smooth map.*

*Define $K_f \subseteq M \times M$ to be the submanifold $\{(x,y) : f(x) = f(y)\}$. We regard this as a space over $M$ via projection onto the first coordinate. We note also that the inclusion of the diagonal in $M \times M$ factors through $K_f$.*

*Define $G_f \subseteq M \times N$ to be the submanifold $\{(x, f(x))\}$. We regard this as a space over $M$ via projection onto the first coordinate.*

*Then there are smoothly locally kinematic structures on $\Delta(M) \subseteq K_f$, $M \subseteq M \times M$, and $G_f \subseteq M \times N$ such that the following diagram commutes.*

$$\begin{array}{ccccc}
M & \longrightarrow & M & \xrightarrow{(1 \times f)\Delta} & G_f \\
\downarrow{\Delta} & & \downarrow{\Delta} & & \downarrow \\
K_f & \longrightarrow & M \times M & \xrightarrow{1 \times f} & M \times N \\
\uparrow & & \uparrow & & \uparrow \\
\ker Tf & \longrightarrow & TM & \xrightarrow{Tf} & f^*TN
\end{array}$$

*Proof.* The key here is to construct a smoothly locally kinematic structure on $M$ taking into account $f$ right from the start. A smoothly locally kinematic structure on $M$ takes tangent vectors to points on the manifold. We want to do this in such a way that the tangent vectors in the same coset of $\ker Tf$ are taken to the same level set of $f$.

We need to choose some additional structure. We need complements of $\ker Tf$ in $TM$ and of $\operatorname{im} Tf$ in $f^*TN$. Let us write these as $E_f$ and $F_f$ respectively. We identify $\operatorname{im} Tf$ with $E_f$ via $Tf$ and so consider the decomposition of $f^*TN$ to be into $E_f \oplus F_f$. We also need to choose connections on $\ker Tf$ and on $F_f$. The connections give rise to parallel transport operators and it is these that we shall use. We also need to assume the existence of a smoothly locally kinematic structure on $N$. We don't need this until we have shown how to construct one on $M$ so we can either apply the same construction to $N$ or we can get one "on the cheap" from a tubular neighbourhood of the diagonal in $N \times N$. (The tubular method is not available to us for $M$ since we need that one to have particular properties.)

Let us begin.

Let $k$ be the rank of $f$. Let $p \in M$ (note that $p$ will play a very minimal role in what follows). The rank theorem for regular maps provides charts $\varphi \colon \mathbb{R}^m \supseteq W \to U \subseteq M$ with $0 \mapsto p$ and $\psi \colon \mathbb{R}^n \supseteq X \to V \subseteq N$ with $0 \mapsto f(p)$, with the property that $\psi^{-1} \circ f \circ \varphi$ is defined and is the projection $(x_1, \ldots, x_m) \mapsto (x_1, \ldots, x_k, 0, \ldots, 0)$. Let us adjust the domains of the charts slightly. Choose an open set of the form $\varepsilon_1 D^k \times \varepsilon_2 D^{m-k}$ with closure inside $W$. By construction, the closure of $\varepsilon_1 D^k \times \{0\}$ is in $X$. Choose $\varepsilon_3 > 0$ such that $\varepsilon_1 D^k \times \varepsilon_3 D^{n-k}$ is inside $X$. Now choose diffeomorphisms $\mathbb{R}^k \to \varepsilon_1 D^k$, $\mathbb{R}^{m-k} \to \varepsilon_2 D^{m-k}$, and $\mathbb{R}^{n-k} \to \varepsilon_3 D^{n-k}$. These fit together to give diffeomorphisms $\mathbb{R}^m \to \varepsilon_1 D^k \times \varepsilon_2 D^{m-k}$ and $\mathbb{R}^n \to \varepsilon_1 D^k \times \varepsilon_3 D^{n-k}$. By pre- and post-composing the original charts with these diffeomorphisms (or inverse as appropriate), we get new charts near $p$ and $f(p)$ with domains $\mathbb{R}^m$ and $\mathbb{R}^n$ such that $f$ is again the projection onto the first $k$ coordinates.

To summarise, we have charts $\varphi_p \colon \mathbb{R}^m \to U_p \subseteq M$ and $\psi_p \colon \mathbb{R}^n \to V_p \subseteq M$ with $\varphi_p(0) = p$ and $\psi_p(0) = f(p)$, such that $\psi_p^{-1} \circ f \circ \varphi_p(x_1, \dots, x_n) = (x_1, \dots, x_k, 0, \dots, 0)$.

Using the chart for $M$ we define a smooth map $X_p \colon TU_p \to \mathcal{X}(U_p)$ as follows. Under the chart map this corresponds to a smooth map $T\mathbb{R}^m \to \mathcal{X}(\mathbb{R}^m)$. Furthermore, we identify $T\mathbb{R}^m$ with $\mathbb{R}^m \times \mathbb{R}^m$ and $\mathcal{X}(\mathbb{R}^m)$ with $C^\infty(\mathbb{R}^m, \mathbb{R}^m)$. We want the resulting map $\mathbb{R}^m \times \mathbb{R}^m \to C^\infty(\mathbb{R}^m, \mathbb{R}^m)$ to be $(q, v) \mapsto (x \mapsto v)$. The key properties that $X_p$ has are:

1. For $v \in TU_p$, $X_p(v)(\pi(v)) = v$.

2. If $Tf(v) = 0$ then $Tf(X_p(v)(x)) = 0$.

By construction the family $\{U_p\}$ covers $M$. We choose a variant of a partition of unity subordinate to this family. The variation is that we want the sum of the squares of the functions to sum to 1. This is a simple modification of the standard construction. Thus we choose a family of functions $\{\rho_i \colon M \to \mathbb{R}\}_{i \in I}$ with the following properties:

1. The closures of the supports of the $\rho_i$ form a locally finite cover of $M$ which refines $\{U_p\}$.

2. The closures of the supports of the $\rho_i$ are compact.

3. For each $x \in M$, $\sum \rho_i(x)^2 = 1$.

For $i \in I$ choose $p_i \in M$ with the closure of the support of $\rho_i$ contained in $U_{p_i}$. Let $U_i := U_{p_i}$. Let $\pi \colon TM \to M$ be the projection. Define $X_i \colon TM \to \mathcal{X}(M)$ by the formula:

$$X_i(v)(q) = \begin{cases} \rho_i(\pi(v))\rho_i(q)X_{p_i}(v)(q) & \text{if } \pi(v), q \in U_i \\ 0 & \text{otherwise} \end{cases}$$

This is smooth as $X_i(v) = 0$ when $\rho_i(\pi(v)) = 0$. It is also the case that $X_i(v)$ has compact support since $X_i(v)(q) = 0$ when $\rho_i(q) = 0$.

As the support of $X_i$ is contained in that of $\rho_i$ we can define a smooth function $X \colon TM \to \mathcal{X}(M)$ by $X(v) = \sum X_i(v)$. Since each $X_i(v)$ has compact support and, for a given $v$, only a finite number of them are non-zero then $X(v)$ has compact support. Thus the function actually has codomain in $\mathcal{X}_c(M)$.

This function inherits the following properties from the $X_p$:

1. For $v \in TM$,

$$X(v)(\pi(v)) = \sum \rho_i(\pi(v))\rho_i(\pi(v))X_{p_i}(v)(\pi(v)) = \sum \rho_i(\pi(v))^2 v = v.$$

2. For $v \in \ker Tf$,

$$Tf(X(v)(x)) = \sum \rho_i(\pi(v))\rho_i(x)TfX_{p_i}(v)(x) = \sum \rho_i(\pi(v))\rho_i(x)0 = 0.$$

Let exp: $\mathfrak{X}_c(M) \to \text{Diff}(M)$ denote the exponential function. Let $\eta_1 : TM \to M$ be the function $\eta_1(v) = \exp(X(v))(\pi(v))$. That is, we take the integral curve of $X(v)$ with initial position at $\pi(v)$ and then evaluate at time 1.

Now if $v \in \ker Tf$ then for $x \in M$, $TfX(v)(x) = 0$. Thus $f$ is constant on the integral curves of $X(v)$. In particular, $f(\eta_1(v)) = f(\pi(v))$. Hence when restricted to $\ker Tf$, $\pi \times \eta_1$ takes values in $K_f$.

We want to modify $\eta_1$ slightly. At the moment, $\pi \times \eta_1$ takes $\ker Tf$ to $K_f$ but the modification will mean that it also takes all cosets of $\ker Tf$ to $K_f$. To do this, we need the decomposition $TM \cong \ker Tf \oplus E_f$ and the parallel transport operator on $\ker Tf$. For a path $\alpha : [0,1] \to M$ let us write parallel transport along $\alpha$ as $\|_\alpha$. We define $\eta_2 : TM \to M$ as the following composition:

1. Start with $v \in TM$.

2. We apply the decomposition $TM \cong \ker Tf \oplus E_f$ to get $v_K \in \ker Tf$ and $v_E \in E_f$.

3. For $t \in [0,1]$, we apply $\eta_1$ to $tv_E$ to get a path, say $\alpha$, in $M$ from $\pi(v)$ to $\eta_1(v_E)$.

4. We parallel transport $v_K$ along this path to $\|_\alpha v_K$ in $\ker Tf$ above $\eta_1(v_E)$.

5. We apply $\eta_1$ to $\|_\alpha v_K$.

The point of this construction is that $f\eta_2(v_K + v_E) = f\eta_2(v_E)$ so the cosets of $\ker Tf$ get mapped to the level sets of $f$.

At this point, we turn our attention to $N$. Here we need the decomposition $f^*TN \cong E_f \oplus F_f$ and the parallel transport operator on $F_f$. We shall use the same notation for this. We also need the smoothly locally kinematic structure on $N$. This defines a smooth map $TN \to N \times N$. Let $\eta_3 : TN \to N$ be the composition of this with the projection onto the second coordinate.

We define $\eta_4 : f^*TN \to N$ as the following composition.

1. Start with $v \in f^*TN$.

2. We apply the decomposition $f^*TN \cong E_f \oplus F_f$ to get $v_E \in E_f$ and $v_F \in F_f$.

3. For $t \in [0,1]$, we apply $\eta_1$ to $tv_E$ to get a path, say $\alpha$, in $M$ from $\pi(v)$ to $\eta_1(v_E)$.

4. We parallel transport $v_F$ along this path to $\|_\alpha v_F$.

5. We push this forward to the vector $Tf \|_\alpha v_F$ in $TN$ above $f(\eta_1(v_E))$.

6. We apply $\eta_4$ to $Tf \|_\alpha v_F$.

The point of this construction is that $f\eta_2(v_E + v_K) = \eta_4(Tfv_E)$.

Let $\nu_M : TM \to M \times M$ be $\pi \times \eta_2$. By construction $\nu_M$ restricts to a map $\ker Tf \to K_f$ which we shall denote by $\nu_K$. Let $\nu_N : f^*TN \to M \times N$ be $\pi \times \eta_4$.

Let $\Delta : M \to M \times M$ be the diagonal map. The image of this map lies in $K_f$ and we write $\Delta$ for the map $M \to K_f$ as well. Then $\Delta^*T(M \times M) \cong TM \oplus TM$ and $\Delta^*TK_f \cong TM \times \ker Tf$. Let $\Gamma : M \to M \times N$ be the graph map, $x \mapsto (x, f(x))$. Then $\Gamma^*T(M \times N) \cong TM \oplus f^*TN$. Let $\zeta : M \to TM$ be the zero section, and we use $\zeta$ also for the zero sections of $\ker Tf$ and $f^*TN$. Then $\zeta^*T(TM) \cong TM \oplus$ $TM$ where the factors are horizontal and vertical. Also $\zeta^*T(\ker Tf) \cong TM \oplus \ker Tf$ and $\zeta^*T(f^*TN) \cong TM \oplus f^*TN$. The maps $\nu_K$, $\nu_M$, and $\nu_N$ have all been constructed so that they induce the obvious isomorphisms $\zeta^*T(TM) \cong \Delta^*T(M \times M)$, $\zeta^*T(\ker Tf) \cong \Delta^*TK_f$, and $\zeta^*T(f^*TN) \cong \Gamma^*T(M \times N)$. Hence they restrict to diffeomorphisms on neighbourhoods of the zero sections.

The final step is to scale the bundles into the neighbourhoods where these maps are diffeomorphisms. We want to do this compatibly on all the bundles involved so we shall scale the bundles $\ker Tf$, $E_f$, and $F_f$ and then scale $TM$ and $f^*TN$ by projecting to the appropriate subbundles first.

The result is maps $\ker Tf \to K_f$, $TM \to M \times M$, and $f^*TN \to M \times N$ which are diffeomorphisms onto their images, being open neighbourhoods of the diagonal (or graph), and which are compatible with $f$ in the appropriate sense.

**Corollary 5.2** *Let $f \colon M \to N$ be a regular map between finite dimensional smooth manifolds. Let $S$ be a smoothly $\mathfrak{T}$–compact space. Then $C^\infty(S, f) \colon C^\infty(S, M) \to C^\infty(S, N)$ is a regular map.*

Buried in the proof of Theorem 5.1 is the result that bundles with connection are smoothly locally trivial. Let us extract that.

**Proposition 5.3** *Let $M$ be a smooth manifold. Let $P \to M$ be a fibre bundle which admits a parallel transport operator. Then $P$ is smoothly locally trivial over $M$.*

*Proof.* Let $\nu \colon TM \to V \subseteq M \times M$ be a smoothly locally kinematic structure on $M$. Then scalar multiplication on $TM$ defines for $(x, y) \in V$ a path $\alpha \colon [0, 1] \to M$ with $\alpha(0) = x$ and $\alpha(1) = y$. Let $\|_{(x,y)} \colon P_x \to P_y$ be the isomorphism coming from parallel transport along this path. Then $(x, y) \mapsto \|_{(x,y)}$ defines a bundle isomorphism $\pi_1^* P \cong \pi_2^* P$ as required.

To apply the work of Section 3.2, we need to know that submanifolds are smoothly locally deformable in their ambient spaces.

**Proposition 5.4** *Let $M$ be a smooth finite dimensional manifold and let $\pi \colon E \to M$ be a finite dimensional vector bundle over $M$. Then there is a smooth map $\psi \colon \mathbb{R} \times E \to \mathrm{Diff}_c(E)$ (diffeomorphisms with compact support) with the following properties:*

1. *$t \mapsto \psi(t, v)$ is a group homomorphism (with $\mathbb{R}$ additive),*

2. *$\psi(t, v)$ is fibre-preserving; that is, for all $x \in E$, $\pi\psi(t, v)(x) = \pi(x)$,*

3. *Let $\zeta \colon M \to E$ denote the zero section. Then $\psi(1, v)(v) = \zeta\pi(v)$.*

*Proof.* As we have already used the exponential map from vector fields to diffeomorphisms elsewhere we shall use it again for this result. For a more direct construction see [nLa11]. The key to this is to construct suitable vector fields on the total space of $E$ and then apply the exponential map to diffeomorphisms.

As $E$ is a finite dimensional vector bundle over $M$ it is smoothly locally trivial. Let $V \subseteq M \times M$ be a suitable neighbourhood of the diagonal and let $\Psi \colon \pi_1^* E \to \pi_2^* E$ be a bundle isomorphism over $V$. The topology on $M \times M$ is smoothly regular and so there is a smooth function $\sigma \colon M \times M \to [0, 1]$ which is 1 on the diagonal and such that the closure of the support of $\sigma$ is contained in $V$.

Using these we construct a smooth map $\overline{X} \colon E \to \mathcal{X}(E)$ as follows. Let $v \in E$. Let $u \in E$ with $(\pi v, \pi u) \in V$. By assumption $\Psi_{(\pi v, \pi u)}$ is an isomorphism $E_{\pi v} \to E_{\pi u}$ so we can use this to transfer $v$ to $E_{\pi u}$. Then we use the cannonical map $\iota\colon \pi^*E \to TE$ to transfer this to $T_uE$ where it is a vertical tangent vector based at $u$. To extend this to all of $E$ we multiply this by $\sigma(\pi v, \pi u)$. Thus $\overline{X}$ is defined by:

$$\overline{X}(v)(u) = \begin{cases} \sigma(\pi v, \pi u)\iota\Psi_{(\pi v, \pi u)}(v) & \text{if } (\pi v, \pi u) \in V \\ 0 & \text{otherwise.} \end{cases}$$

This is not quite the map we want as we want it to be compactly supported also in the fibre directions. The key property of this is that for $t \in [0,1]$ then $\overline{X}(v)(tv) = v$. We want to keep this, which means that we can only modify $\overline{X}(v)$ away from $v$. A simple way to do this is to choose a smooth orthogonal structure on $E$ and a function $\tau\colon \mathbb{R} \to [0,1]$ such that for $t \in [0,1]$ then $\tau(t) = 1$ whilst for $t \geq 2$ then $\tau(t) = 0$. Let $\rho\colon E \to \mathbb{R}$ be the quadratic form associated to the orthogonal structure. Then we define $X\colon E \to \mathcal{X}(E)$ by:

$$X(v)(u) = \tau\left(\frac{\rho(u)}{1+\rho(v)}\right)\overline{X}(v)(u).$$

At a vector $v$ the resulting vector field is then compactly supported, both on the base and in the fibres. For $t \in [0,1]$ then $X(v)(tv) = \iota v$.

We now apply the exponential map to produce a smooth map $\psi\colon \mathbb{R} \times E \to \text{Diff}_c(E)$, $\psi(t,v) = \exp(-tX(v))$ (note the minus sign). By construction for $v \in E$ then the map $t \mapsto \psi(t,v)$ is a group homomorphism. As $X(v)$ is a vertical vector field $\pi\psi(t,v)$ is independent of $t$ and thus $\pi\psi(t,v) = \pi$.

Lastly, let us consider $\psi(1,v)(v)$. To determine this we need to look at the integral curve of $-X(v)$ starting at $v$. As $X(v)$ is a vertical vector field we only need to look in the fibre of $E_{\pi v}$ and, at least to start with, we are within the region where $X(v)(u) = \iota v$. The integral curve is thus $t \mapsto v - tv$. Whence $\psi(1,v)(v) = \zeta\pi(v)$. 🌿

**Corollary 5.5** 1. *Let $N$ be a smooth finite dimensional manifold and let $M \subseteq N$ be a submanifold. Then $M$ is smoothly locally deformable in $N$.*

2. *Let $M$ be a smooth finite dimensional manifold. Then $M$ is smoothly locally deformable.* 🌿

This means that we can apply the results of Section 3.2 to the case of smooth manifolds.

**Theorem 5.6** *Let $N$ be a smooth finite dimensional manifold and $M \subseteq N$ an embedded submanifold. Let $E \to M$ be the normal bundle of the embedding. Let $S$ be a compact smooth manifold and $T \subseteq S$ an embedded submanifold. Then $C^\infty((T,S),(M,N))$ is a submanifold of $C^\infty(S,N)$ which admits a tubular neighbourhood with normal bundle the pull-back of $C^\infty(T,E)$.*

*Proof.* The requirements of Proposition 3.10 are satisfied. 🌿

Using this we recover the following result of string topology.

**Corollary 5.7** *Let $M$ be a smooth manifold. Let $8_\times$ denote the figure eight formed by taking the wedge of two circles. Then $C^\infty(8_\times, M)$ is a submanifold of $C^\infty(S^1, M) \times C^\infty(S^1, M)$ with a tubular neighbourhood diffeomorphic to the pull-back of $TM$ under the evaluation map at the basepoint $C^\infty(8_\times, M) \to M$.* 🌿

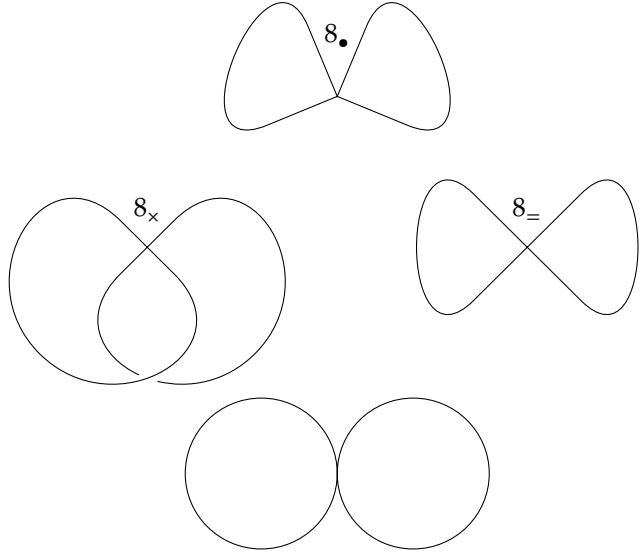

Figure 1: Four smooth structures on the figure 8

It is obvious how to extend this to a cactus.

Using quotients instead of subobjects we deduce the following result.

**Theorem 5.8** *Let $M$ be a smooth manifold. Let $T$ be a compact smooth manifold. Let $R \subseteq T$ be an embedded submanifold. Let $\rho\colon R \to T$ be a smooth map such that the closure of the image of $R$ does not meet the closure of $R$. Let $S$ be the coequaliser of $R \rightrightarrows T$. Then $C^\infty(S, X)$ is a submanifold of $C^\infty(T, X)$ which admits a tubular neighbourhood with normal bundle the pull-back of $C^\infty(R, TM)$.*

*Proof.* The requirements of Proposition 3.12 are satisfied.

This allows us to recover the opposite result of string topology.

**Corollary 5.9** *Let $M$ be a smooth manifold. Let $8_=$ denote the figure eight formed by identifying two antipolar points on a circle. Then $C^\infty(8_=, M)$ is a submanifold of $C^\infty(S^1, M)$ with a tubular neighbourhood diffeomorphic to the pull-back of $TM$ under the evaluation map at the identified point $C^\infty(8_=, M) \to M$.*

Note that the two spaces $C^\infty(8_=, M)$ and $C^\infty(8_\times, M)$ are not directly related. They both embed as submanifolds of the manifold $C^\infty(8_\bullet, M)$ where $8_\bullet$ is the figure eight formed by taking two intervals and identifying all of the endpoints. Their intersection is a little trickier to describe, essentially the two branches of the space can be chosen independently at the singular point. See Figure 1 for an approximation of the four structures. Each submanifold is of infinite codimension in its containing manifold. Thus to pass from $C^\infty(8_=, M)$ to $C^\infty(8_\times, M)$ we have to pass up and down an infinite number of dimensions.

Lastly, let us observe that any smoothly locally kinematic structure on a finite dimensional manifold is Fermat. (Recall that we've assumed that the kinematic tangent space of a smooth manifold is its standard tangent space.)

**Proposition 5.10** *Let M be a finite dimensional smooth manifold. Then a smoothly local deflation on M which is kinematic is also Fermat.*

*Proof.* As the smoothly local deflation is kinematic we can transfer the target of the Fermat map to the tangent bundle of *M*. Since we've assumed that the kinematic tangent functor agrees with the standard one locally this looks like $\mathbb{R}^n \times \mathbb{R}^n$. We can therefore use Hadamard's construction from Equation (4) to define the function explicitly. 🐦

## 6 Frölicher Spaces

In this section we shall consider how some of the results of this article look in the category of Frölicher spaces, which is a suitable cartesian closed category of generalised smooth spaces.

**Definition 6.1** *A* Frölicher space *consists of a triple* $(X, C, \mathcal{F})$ *where X is a set,* $C \subseteq \text{Set}(\mathbb{R}, X)$ *a family of curves, and* $\mathcal{F} \subseteq \text{Set}(X, \mathbb{R})$ *a family of functionals. These satisfy the following equivalent:*

1. *A curve* $c \colon \mathbb{R} \to X$ *is in* $C$ *if and only if for all* $f \in \mathcal{F}$ *then* $f \circ c \in C^\infty(\mathbb{R}, \mathbb{R})$*, and*

2. *A functional* $f \colon X \to \mathbb{R}$ *is in* $\mathcal{F}$ *if and only if for all* $c \in C$ *then* $f \circ c \in C^\infty(\mathbb{R}, \mathbb{R})$.

*The elements of* $C$ *are called* smooth curves *in X, whilst of* $\mathcal{F}$ *are* smooth functionals *on X.*

*A morphism of Frölicher spaces is a set function* $\varphi \colon X \to Y$ *which satisfies one (hence all) of the following equivalent conditions:*

1. *for all* $c \in C_X$ *then* $\varphi \circ c \in C_Y$,

2. *for all* $f \in \mathcal{F}_Y$ *then* $f \circ \varphi \in \mathcal{F}_X$,

3. *for all* $c \in C_X$ *and* $f \in \mathcal{F}_Y$ *then* $f \circ \varphi \circ c \in C^\infty(\mathbb{R}, \mathbb{R})$.

As with any category of generalised smooth spaces we have the two extreme topologies: the *curved* and *functional* ones. For the curved topology we have the following characterisation of the smoothly compact spaces.

**Proposition 6.2** *A Hausdorff Frölicher space is smoothly compact for the curved topology if and only if its curved topology is sequentially compact.*

The difference between a Hausdorff Frölicher space and an arbitrary one is not a lot. Every Frölicher space has a Hausdorffification and the fibres of the projection are equipped with the indiscrete topology.

*Proof.* With the curved topology a subset is open if and only if its preimage under all smooth curves is open in $\mathbb{R}$. So let *S* and *X* be arbitrary Frölicher spaces. Let $U \subseteq X$ be a subset that is curve open. We wish to determine when $C^\infty(S, U)$ is curve open in $C^\infty(S, X)$.

Let $c \colon \mathbb{R} \to C^\infty(S, X)$ be a smooth curve. Suppose, without loss of generality, that $c(0) \in C^\infty(S, U)$. To show that $C^\infty(S, U)$ is curve open in $C^\infty(S, X)$ we need to find some $\varepsilon > 0$ such that for $|t| < \varepsilon$ then $c(t) \in C^\infty(S, U)$.

As the category of Frölicher spaces is cartesian closed the curve *c* defines a smooth map $\check{c} \colon \mathbb{R} \times S \to X$ by $\check{c}(t, s) = c(t)(s)$. To say that for $|t| < \varepsilon$ that $c(t) \in C^\infty(S, U)$ means that for $(t, s) \in (-\varepsilon, \varepsilon) \times S$ then $\check{c}(t, s) \in U$. Thus $\varepsilon$ exists if and only if $\check{c}^{-1}(U)$ contains a subset of the form $(-\varepsilon, \varepsilon) \times S$.

As assigning the curved topology is functorial and $U$ is curve open in $X$, $\check{c}^{-1}(U)$ is curve open in $\mathbb{R} \times S$. Since we have assumed that $c(0) \in C^\infty(S, U)$ we know that $\{0\} \times S \subseteq \check{c}^{-1}(U)$. Thus we need to find conditions on $S$ that mean that every open neighbourhood of $\{0\} \times S$ contains a subset of the form $(-\varepsilon, \varepsilon) \times S$. Since we are dealing with the curved topology we test open sets by smooth curves.

Let us start by assuming that every open neighbourhood of $\{0\} \times S$ in $\mathbb{R} \times S$ contains a subset of the form $(-\varepsilon, \varepsilon) \times S$. Let $(s_n)$ be a sequence in $S$. Consider the set $W \subseteq \mathbb{R} \times S$ defined to be the complement of the set $\{(1/n, s_n)\}$. By construction this does not contain a subset of the form $(-\varepsilon, \varepsilon) \times S$ and so cannot be an open neighbourhood of $\{0\} \times S$. It must therefore be the case that there is a smooth curve $c \colon \mathbb{R} \to \mathbb{R} \times S$ such that $c^{-1}(W)$ is not an open subset of $\mathbb{R}$. We can therefore find a sequence $(t_k)$ in $\mathbb{R}$ with the property that for each $k$, $t_k \notin c^{-1}(W)$ but $(t_k)$ converges to a point, say $t$, in $c^{-1}(W)$. Since $t_k \notin c^{-1}(W)$ we must have that there is some $n_k$ such that $c(t_k) = (1/n_k, s_{n_k})$. Since $(t_k)$ converges to a point in $c^{-1}(W)$ the set $\{n_k\}$ must be infinite, whence unbounded. We can therefore find a subsequence, say $(t_{k_j})$, such that the corresponding sequence $(n_{k_j})$ is strictly increasing. The sequence $(t_{k_j})$ still converges to $t$ and so $(c(t_{k_j}))$ converges to $c(t)$. Its projection to $S$ therefore converges in $S$. This projection is $(s_{n_{k_j}})$ and therefore $(s_n)$ has a convergent subsequence. Hence $S$ is sequentially compact.

Now consider the converse. For this we need to assume that $S$ is Hausdorff. Suppose that $S$ is not sequentially compact. Then there is a sequence, say $(s_n)$, which does not have a convergent subsequence. As above, let $W$ be the complement in $\mathbb{R} \times S$ of the set $\{(1/n, s_n)\}$. We wish to show that this is open in $\mathbb{R} \times S$. So let $c \colon \mathbb{R} \to \mathbb{R} \times S$ be a smooth curve and consider $c^{-1}(W)$. The complement of this is the set $\{t \in \mathbb{R} \colon c(t) \in \{s_n\}\}$. As $S$ is Hausdorff for each $n$ then the set $c^{-1}(\{s_n\})$ is closed in $\mathbb{R}$. A compact subset of $\mathbb{R}$, say $C$, can only meet a finite number of the sets $c^{-1}(\{s_n\})$ since otherwise there would be a convergent subsequence. Hence the intersection $C \cap c^{-1}(W)$ is open in $C$. Thus as $\mathbb{R}$ is compactly generated $c^{-1}(W)$ is open in $\mathbb{R}$. Thus $W$ is open in $\mathbb{R} \times S$ but does not contain a subset of the form $(-\varepsilon, \varepsilon) \times S$. ◻

Another part of the theory developed here that has some loose ends is that of tangent spaces. The following result on tangent spaces of functionally open subsets is useful in that regard.

**Proposition 6.3** *Let $X$ be a Frölicher space. Let $U \subseteq X$ be a functionally open subset. Then $TU \to TX$ is an embedding (onto a functionally open subset).*

*Proof.* That the target is functionally open is immediate as it is the preimage of $U$ under the projection $TX \to X$. The map $TU \to TX$ is smooth and is a bijection by Proposition 4.7.

What remains is to show that it is an embedding. To do this we need to show that if $\alpha \colon \mathbb{R} \to TX$ is a smooth curve with the property that its image is in $TU$ then it is smooth into $TU$. To prove this we need to examine functionals on $TU$.

Let $f \colon TU \to \mathbb{R}$ be smooth. Let $v \in TU$. Then $\pi v \in U$ and as $U$ is functionally open in $X$ there is some smooth $g \colon X \to \mathbb{R}$ with the property that $g(x) = 1$ and $g$ is zero outside $U$. Let $h$ be the composition of $g$ with a suitable smooth function so that $h$ takes values in $[0,1]$ and is $1$ on the set $\{y \colon g(y) > {}^3\!/\!_4\}$ and is $0$ on the set $\{y \colon g(y) < {}^1\!/\!_4\}$ (as in the proof of Proposition 4.7). Consider the product $f \bullet (h \circ \pi)$, where $\bullet$ stands for pointwise multiplication of functions. This is a smooth function $TU \to \mathbb{R}$ which agrees with $f$ on the set $\{u \colon g \circ \pi(u) > {}^3\!/\!_4\}$, which is a neighbourhood of $v$, and is zero on the set $\{u \colon g \circ \pi(u) < {}^1\!/\!_4\}$.

Now we define a function $\hat{f} \colon TX \to \mathbb{R}$ by defining it to be $f \bullet (h \circ \pi)$ on $TU$ and $0$ elsewhere. We wish to show that this is smooth. As $TX$ is a quotient of $C^\infty(\mathbb{R}, X)$, $\hat{f}$ is smooth on $TX$ if and only if its composition with the projection is smooth. Let $q \colon C^\infty(\mathbb{R}, X) \to TX$ be the quotient. To show that $\hat{f} \circ q$ is smooth we consider a smooth curve $\beta \colon \mathbb{R} \to C^\infty(\mathbb{R}, X)$ and the composition $\hat{f} \circ q \circ \beta \colon \mathbb{R} \to \mathbb{R}$.

To show that this is smooth it is sufficient to work locally. So let $t \in \mathbb{R}$ and consider $\beta(t)$. This is a curve in $X$. Suppose that $\beta(t)(0) \notin U$, whence $g(\beta(t)(0)) = 0$. The map $C^\infty(\mathbb{R}, X) \to X$ given by evaluation at $0$ is smooth, whence continuous, so there is some $\varepsilon > 0$ such that $g(\beta(s)(0)) < {}^1\!/\!_4$ for $|s - t| < \varepsilon$. Then $\hat{f} \circ q \circ \beta$ is $0$ in a neighbourhood of $t$ and thus is smooth near $t$.

Now suppose that $\beta(t)(0) \in U$. If $g(\beta(t)(0)) < {}^1\!/\!_4$ then we can argue as above. So assume that $g(\beta(t)(0)) \geq {}^1\!/\!_4$. As we're working in a cartesian closed category $\beta$ defines a smooth map $\hat{\beta} \colon \mathbb{R}^2 \to X$. Since $\hat{\beta}(t, 0) = \beta(t)(0) \in U$ there is some $\varepsilon > 0$ such that if $|s - t| < \varepsilon$ and $|r| < \varepsilon$ then $\hat{\beta}(s, r) \in U$ (the functional topology on $\mathbb{R}^2$ agrees with the standard one).

Let $\nu \colon \mathbb{R} \to \mathbb{R}$ be a smooth map with the property that for all $r \in \mathbb{R}$ then $|\nu(r)| < \varepsilon$ whilst for $|r| \leq \varepsilon/2$ then $\nu(r) = r$. Let $\tilde{\beta} \colon \mathbb{R} \to C^\infty(\mathbb{R}, X)$ be the function $\tilde{\beta}(s)(r) = \beta(\nu(s))(\nu(r))$. By construction for all $s, r \in \mathbb{R}$ then $\tilde{\beta}(s)(r) \in U$ whilst for $s, r$ with $|s - t|, |r| \leq \varepsilon/2$ then $\tilde{\beta}(s)(r) = \beta(s)(r)$.

Since for $s \in \mathbb{R}$, $\tilde{\beta}(s)$ takes values in $U$ it is actually a smooth curve in $U$ whence $\tilde{\beta}$ can be considered as a smooth map $\mathbb{R} \to C^\infty(\mathbb{R}, U)$ (as the associated map $\mathbb{R}^2 \to U$ is smooth). We therefore have that $(f \bullet (h \circ \pi)) \circ q \circ \tilde{\beta}$ is a smooth map $\mathbb{R} \to \mathbb{R}$. For $|s - t| < \varepsilon/2$ and $|r| < \varepsilon/2$ then $\tilde{\beta}(s)(r) = \beta(s)(r)$, whence for $|s - t| < \varepsilon/2$ then $q \circ \tilde{\beta}(s) = q \circ \beta(s)$. Hence $\hat{f} \circ q \circ \tilde{\beta}$ agrees with $(f \bullet (h \circ \pi)) \circ q \circ \beta$ in a neighbourhood of $t$ and thus is smooth in that neighbourhood.

Thus $\hat{f} \circ q \circ \beta$ is smooth, whence $\hat{f}$ is smooth.

Now for $\alpha \colon \mathbb{R} \to TX$ with the property that its image is in $TU$ consider some $t \in \mathbb{R}$ and let $v = \alpha(t)$. Let $f \colon TU \to \mathbb{R}$ be smooth. The construction above yields a smooth function $\hat{f} \colon TX \to \mathbb{R}$ which agrees with $f$ on a neighbourhood of $v$. Hence $f \circ \alpha$ is smooth in a neighbourhood of $t$ and thus is a smooth function. Therefore $\alpha$ is a smooth curve in $TU$ and so $TU \to TX$ is an embedding. 🐦

This has important consequences for smoothly locally kinematic spaces.

**Theorem 6.4** *Let X be a smoothly locally kinematic space with respect to the functional topology functor. Then X is smoothly locally additive and (using the notation established in Section 4) there is a diffeomorphism $\iota^* \ker T\tau \cong TX$.*

*Proof.* We assume the notation of Proposition 4.12 and its proof. Since $T(\mathbb{R}^2 \times V \times_X V) \cong T\mathbb{R}^2 \times T(V \times_X V)$ we can define a smooth map $\mathbb{R}^2 \times (\mathbb{R}^2 \times V \times_X V) \to T(\mathbb{R}^2 \times V \times_X V)$ of spaces over $\mathbb{R}^2 \times V \times_X V$ where the $\mathbb{R}^2$–factor is the $\mathbb{R}^2$–tangent space. As $W \subseteq \mathbb{R}^2 \times V \times_X V$ is open $TW$ is an embedded subspace of $T(\mathbb{R}^2 \times V \times_X V)$. Hence the restriction of this section defines a smooth map $\mathbb{R}^2 \times W \to TW$. This is what is needed to make the addition of vectors smooth.

To show that $\iota^* \ker T\tau \to TX$ is a diffeomorphism we observe first that $TV \to T(X \times X) \cong TX \times TX$ is an embedding. This identifies $\ker T\tau$ with the second factor, whence $\iota^* \ker T\tau$ is diffeomorphic to $TX$.

The other aspect of the theory that has been developed here which is somewhat open is that of when a mapping space is smoothly locally kinematic.

**Theorem 6.5** *Let M be a finite dimensional smooth manifold. Let S be a smoothly $\mathfrak{T}_f$–compact space. Then $C^\infty(S, M)$ is smoothly locally kinematic.*

*Proof.* Finite dimensional smooth manifolds are Fermat spaces and Theorem 6.4 shows that the technical condition for Theorem 4.15 is satisfied. Also, Theorem 6.4 shows that the conclusion of Theorem 4.15 is sufficient to show that $TC^\infty(S, X) \to C^\infty(S, TX)$ is injective, whence it is a diffeomorphism as both are quotients of $C^\infty(S \times \mathbb{R}, X)$. Hence from Section 4.3.3, $C^\infty(S, M)$ is again smoothly locally kinematic.